\documentclass[reqno,11pt]{amsart}

\title[Tip multifractal
spectrum for SLE]{Almost sure  multifractal spectrum for 
the tip of an SLE curve}
\author[Johansson Viklund]{Fredrik Johansson Viklund}
\address{Fredrik Johansson Viklund\\
Department of Mathematics\\
Columbia University\\
2990 Broadway\\
New York, NY 10027, USA}
\thanks{Johansson Viklund acknowledges the support of the Simons Foundation and the Knut and Alice Wallenberg Foundation.}
\email{fjv@math.columbia.edu}

\author[Lawler]{Gregory F. Lawler}
\address{Gregory F. Lawler\\ 
Department of Mathematics and Department of Statistics\\
University of Chicago\\
5734 S. University Avenue\\
Chicago, IL 60637, USA}
\thanks{Lawler is supported by National Science
Foundation  grant DMS-0907143.}
\email{lawler@math.uchicago.edu}

\usepackage{verbatim}
\usepackage{amssymb,amsmath,amsthm}
\usepackage{graphicx}
\usepackage{epsfig}
\newtheorem*{theorem*}{Theorem}
\newtheorem{corollary}{Corollary}[section]
\newtheorem{lemma}[corollary]{Lemma}
\newtheorem{proposition}[corollary]{Proposition}

\newtheorem{theorem}[corollary]{Theorem}
\newcommand{\SLE}{\rm{SLE}}

\newcommand{\Prob} {{\mathbb P}}
\newcommand{\Z}{{\mathbb Z}} 
\newcommand{\HH}{{\mathbb H}} 
\newcommand{\E}{{\mathbb E}}

\newcommand{\R}{{\mathbb{R}}}
\newcommand{\C}{{\mathbb C}}

\newcommand{\hdim}{{\rm dim}_h}

\newcommand{\dist}{{\rm dist}}

\renewcommand{\Im}{{\rm Im}}
\renewcommand{\Re}{{\rm Re}}
\newcommand{ \p} {\partial}

\newcommand{ \Half }{{\mathbb H}}

\newcommand{ \Disk }{{\mathbb D}}

\newcommand{ \diam }{{\rm diam}}

\newcommand{ \gap} {\succcurlyeq}
\newcommand{ \lap} {\preccurlyeq}
\newcommand{ \gapp} {\succcurlyeq_{\rm i.o.}}
\newcommand{ \lapp}{\preccurlyeq_{\rm i.o.}}
\newcommand{ \wgap}{\succapprox}

\newcommand{ \wgapp}{\wgap_{\rm i.o.}}

\newcommand{ \wapprox}{\approx^*}

\newcommand{ \energy }{{\mathcal E}}

\renewcommand{\k}{{\kappa}}

\newcommand{ \haus }{{\mathcal H}}
\newcommand{ \ball }{{\mathcal B}}
\newcommand{  \supp }{{\rm supp }\,}
\newcommand{\deriv}{\Delta}

\newcommand{\fr}{\frac}
\newcommand{\iy}{\infty}
\newcommand{\ee}{\epsilon}

\newcommand {\hm}{{\rm hm}}
\newcommand{\dyad}{{\mathcal D}}

 \theoremstyle{definition}
\newtheorem*{definition}{Definition}

\theoremstyle{remark}
\newtheorem*{remark}{Remark}
\numberwithin{equation}{section}

\numberwithin{equation}{section}

\begin{document}

\begin{abstract}
The tip multifractal spectrum of a two-dimensional curve
is one way to describe the behavior
of the 
uniformizing conformal map of the complement near the
tip. We give
the tip multifractal spectrum for a Schramm-Loewner
evolution ($\SLE$) curve, we prove that the
spectrum is valid with probability one, and we give applications to the scaling of harmonic measure at the tip.
\end{abstract}

\maketitle

\section{Introduction}

The chordal Schramm-Loewner evolution ($\SLE_\kappa$) is a
one parameter family of  probability
measures on curves $\gamma:[0,\infty) \rightarrow \overline
\Half$, where $\Half$ denotes the complex upper half plane. It was invented by Schramm \cite{Schramm} as
a candidate for the scaling limit of two-dimensional lattice models from statistical physics that satisfy conformal invariance and a Markovian property in
the limit. Several lattice models have since been shown to have scaling limits that can be described by $\SLE$. Examples include loop-erased random walk and the uniform spanning tree \cite{LSW_lerw}, the percolation exploration-process  \cite{Smirnov_perc}, and the FK-Ising model \cite{Smirnov_ising}. We refer the reader to \cite{LConv,LPark,Werner} for surveys and further references. 

In this paper we will be interested in the behavior at
the tip $\gamma(t)$ of the growing $\SLE$ curve.  Since the curves are fractals, one
cannot make sense of derivatives.  Instead, the natural approach is to consider
the behavior of $|g_t'(z)|$ for $z$ near $\gamma(t)$ where
$g_t$ is a uniformizing conformal map from the complement of the
curve to the upper half plane. For technical reasons, it is often easier to consider
$f_t = g_t^{-1}$ near $V_t$, the pre-image of the tip on the real-line. Our main goal will be to derive the almost sure \emph{tip multifractal spectrum} for $\SLE$. For a suitable interval of $\alpha$, it is defined roughly as the dimension of the subset of the curve corresponding to $t$ for which $y|f_t'(iy+V_t)|$ decays like $y^{\alpha}$ when $y \to 0+$. We shall see that the tip multifractal spectrum is closely related to the multifractal spectrum of harmonic measure at the tip. As a function of $\alpha$, this spectrum measures the size of the part of the curve that corresponds to $t$ for which the harmonic measure of a ball of radius $\ee$ centered at the tip decays like $\ee^\alpha$ as $\ee \to 0+$.
 
The multifractal spectrum of harmonic measure has been studied
extensively in the physics and mathematics literature.
For example, in the case of the paths of Brownian motion, the spectrum
is determined by the Brownian intersection exponents, see \cite{LSW1} and the references therein. In two
dimensions these exponents were established by Lawler, Schramm, and Werner in \cite{LSW1,LSW2,LSW3}. In the case of the $\SLE$ path, Duplantier used non-rigorous ``quantum gravity'' methods to predict a harmonic measure spectrum for the tip, see \cite[Section 7]{Duplantier1}. However, this spectrum is different to the ones we we will work with as it describes the \emph{local} dimension of harmonic measure; it corresponds in some sense to our function $\rho(\beta)$, see Section~\ref{SLEsec}. (It also does not consider a ``generic'' tip but rather the non-equivalent behavior at the bulk point of a radial SLE path.) Using similar methods, Duplantier and Binder predicted the spectrum of harmonic measure for the \emph{bulk} of $\SLE$, see \cite{Duplantier_Binder}. Roughly speaking, this spectrum is defined as the dimension of the subset of the curve away from the tip where, for a given $\alpha$, harmonic measure in a shrinking ball of radius $\ee$ decays like $\ee^\alpha$. Beliaev and Smirnov \cite{Beliaev_Smirnov}
made a start to proving this result
by establishing the \emph{average integral means spectrum} for $\SLE$. 
To get the almost sure multifractal spectrum
from the average integral means spectrum one can formally apply the
so-called multifractal formalism \cite{Makarov_fine} and
find the bulk spectrum by taking a Legendre transform of the
average integral means spectrum. This approach is believed to be valid for $\SLE$, 
although it has not been established in this case. Indeed, to the best of our knowledge, our results are the first on almost sure multifractal spectra for the family $\SLE_\kappa, \, \kappa > 0$. 

The starting point of our analysis is estimation
of moments of the derivative of $f_t$ using the reverse-time Loewner flow; this was started
by Rohde and Schramm in \cite{RS} and extended in many places, e.g.,
 \cite{Beliaev_Smirnov,
JL,Kang,Law1,Lind}. (This is the analogue of the
average integral means spectrum result for our problem.) 
In order to get almost sure
results, one needs second moment estimates. The ideas
for that appear in \cite{Law1} and they were used in, e.g., \cite{JL}. These ideas are also
important in understanding the so-called natural parameterization
of $\SLE$ curves, see \cite{LShef}.

\subsection{Multifractal spectra for the tip} We now proceed to discuss in more detail the multifractal spectra that we will consider.
To motivate our definitions, we will start out in a slightly different setting than the one we will work with in the bulk of the paper. 

Suppose that $\zeta$ is a boundary point of a simply
connected domain $D$.  We say that $\zeta$ is {\em accessible (in
$D$ by $\eta$)} if $\eta:[0,1] \rightarrow \C$ is a simple curve 
with $\eta(0) =\zeta$ and $\eta(0,1] \subset D$.  
If $\zeta$ is accessible by $\eta$, let $h$ be a conformal transformation of
$D$ onto $\C \setminus (-\infty,0]$ with $h(\zeta) = 0$.
 By $h(\zeta) = 0$,
we mean $h(\eta(0+)) = 0.$  

We have the following situation in mind. Let $\tilde \gamma:(-\infty, \infty)
\rightarrow \C$ be a simple curve
with $\tilde \gamma(t) \rightarrow \infty$ as $t \rightarrow
\pm \infty$.  For each $t$, we consider the ``slit'' plane $D_t =
\C \setminus \tilde \gamma(-\infty,t]$, which is a simply connected
domain whose boundary contains $\tilde \gamma(t)$ and
$\infty$.  The {\em (nontangential) tip multifractal spectrum} which
we describe in this subsection is
one way to describe 
the behavior near $\tilde \gamma(t)$ of the conformal map uniformizing $D_t$, 
for different values of $t$. Clearly, the boundary point $\tilde \gamma(t)$ is accessible in
$D_t$ by the curve $\eta^{(t)}(s) = \tilde
\gamma(t+s)$.
\begin{remark}For endpoints of slits like $\tilde\gamma(t)$ in
$D_t$, there
is only one possible meaning for $h(\tilde \gamma(t)) = 0$, but 
for general $D$
a boundary point $\zeta$ might be approached from different directions
that correspond to different values of $h(\zeta)$. Formally, this can be understood using \emph{prime ends}, see, e.g., \cite[Chapter 2]{Pommerenke}. In the case at hand, the curve
$\eta$ specifies a particular direction/prime end.
\end{remark}

Returning to the general simply connected domain $D$, let $g(z) = i \, \sqrt {h(z)}$, where the branch of the square
root is chosen so that $\sqrt 1 = 1$. 
Then $g$ is a conformal transformation of $D$ onto the upper half
plane $\Half$ with $g(0) = 0$.  The map $g$ is only unique up to composition with a M\"obius transformation, that is,
 if $\tilde g$ is another such map, then
\[            \tilde g(z) = T[g(z)], \]
where $T$ is a M\"obius transformation of $\Half$
fixing $0$.  Similarly $h$ is not unique.
 
Let $\eta^*(s) = g(\eta(s))$.  Then $\eta^*:(0,1]
\rightarrow \Half$ is a curve with $\eta^*(0+) = 0$.   If $\eta_1^*:(0,1]
\rightarrow \Half$ is another curve with $\eta_1^*(0+) = 0$, 
and $\eta_1(t) = g^{-1}(\eta_1^*(s))$, then 
$\eta_1(0+)= \zeta$ and $\zeta$ is accessible by $\eta_1$. This uses the fact
that the curve $\eta$ exists.  If $\zeta$ is an inaccessible
boundary point, then the limit $ 
\lim_{s \rightarrow 0+}\eta_1(s)$ does not exist. 
We say that $\eta^*$ satisfies a {\em weak cone condition} if
there is a subpower function (see Section \ref{notsec}) $\psi$
such that for all $s > 0$,
\[     \left|\; \Re[\eta^*(s)] \;\right| \leq
     \Im[\eta^*(s)] \, \psi\left(\frac 1{ \Im[\eta^*(s)]}
  \right), \]
and we say that $\eta$ is \emph{weakly nontangential} if
$g \circ \eta$
satisfies a weak cone condition.  It is not difficult to see that this
definition is independent of the choice of $g$. 
 One example of
a weakly nontangential curve for $D$ is
\[     \eta(s) = g^{-1}(si) , \;\;\;\; 0 < s \leq 1.\]
We will use this particular curve to define the tip multifractal spectrum
but the definition will be the same for any weakly nontangential
curve. 

Next, we let $f = g^{-1}$ so that $f$ is a conformal
transformation of $\Half $ onto $D$. 
Since $f(is)=\eta(s), s > 0$, is a simple curve,
  the length of $\eta(0,s]$ is given by
\begin{equation}  \label{nov18.12}
        v(f;s) :=    \int_0^s |f'(iy)| \, dy . 
\end{equation}
A sufficient condition for the existence of a limiting $\zeta
= \eta(0+)$ is
that $v(f;0+) = 0$ which is equivalent to
\[           
v(f;t) < \infty, \quad t > 0. 
\]
We can also use the plane slit by the negative real axis as uniformizing domain and write $f(w) = F(-w^2) $ where $F:
\C \setminus (-\infty,0] \rightarrow D$ with $F(0) = \zeta$.
Then $F^{-1}(\eta(0,s]) = [0,s^2].$
In particular, the length of $F^{-1}(\eta(0,s])$ is $s^2$, and the
length of $f^{-1}(\eta(0,s])$ is $s$.


We will say that  the {\em (nontangential)
scaling exponent} at the boundary point $\zeta$
is $\theta$ if
\[      v(f;s) \wapprox s^{2\theta}, \quad y \rightarrow 0+.\]
In particular, if $D = \C \setminus (-\infty,t]$, then the
scaling exponent at $t$ equals $1$. 
(Recall that $f: \mathbb{H} \to D$ and see \eqref{nov18.1} for the definition of $\wapprox.$)
 More generally, if
$\gamma$ is differentiable at $t$, then $\theta = 1$ at $t$.
Note that the Beurling estimates 
(see Lemma \ref{ber})
imply that $\theta \leq 1$. (In fact, the same bound holds for a ``limsup version'' of the definition of $\theta$.) The scaling exponent is closely
related to the behavior of $|f'(iy)|$ as $ y \rightarrow 0+$.
Indeed, if $y|f'(iy)| \wapprox y^{1-\beta}$
for some $\beta < 1$, then (see Proposition \ref{standard2})
\[   v(f;y) \wapprox y^{1-\beta}, \quad y \to 0+\]
so that
\[
                        \theta = \frac{1-\beta}{2}. \]
Although the definition of $v(f;y)$ depends on the choice of
conformal map $f$, it is not hard to see that the scaling 
exponent $\theta$ is independent of the choice.

Returning to the curve $\tilde \gamma$, we can study $T_\theta$,
the set of $t$
such that the scaling exponent of $D_t$ at $\tilde \gamma(t)$ equals
$\theta$.  
The tip multifractal spectrum can then be defined to be either
of the two functions:
\[ 
\theta \longmapsto \hdim(T_\theta), \;\;\;\;\;
   \theta \longmapsto \hdim(\tilde \gamma[T_\theta]), 
\]
where $\hdim$ denotes Hausdorff dimension.  The first function
depends on the choice of parameterization of $\tilde \gamma$ and the second
is independent of parameterization. One could also define
liminf and limsup versions of this. The main goal of
this paper is to compute
the tip multifractal spectrum for the chordal
$\SLE$ path. For technical convenience, we  will use
 an alternative definition in terms of the behavior of
$|f'(iy)|$ as $y \rightarrow 0+$ and we will use $\beta$ rather
than $\theta$ as our variable.  

Suppose now that $\gamma=\gamma(t)$ is a curve in $
\overline \HH$ with $\gamma(0+)
 \in \R$. 
  Let $H_t$ be the unbounded connected component of $\HH \setminus
\gamma[0,t]$. 
 One way to define 
the multifractal spectrum of harmonic measure at the tip is as
the function
\[
\alpha \longmapsto \hdim(\gamma[T^{\hm}_{\alpha}]),
\]
where
\[
T^{\hm}_{\alpha}=\{t: \, \hm(\mathcal{B}(\gamma(t),r)) \wapprox
r^{\alpha}, \, r \to 0+\}.
\]
Here $\hm(\cdot)=\hm(\infty, \cdot, H_t)$, is the normalized harmonic measure
from infinity. We will 
both use this definition 
and  a slightly
different (nonequivalent)
definition that is more closely related to the tip multifractal spectrum that we described above. See Section \ref{hmsec} for precise definitions.

\subsection{Main results}
The following is part of our main result, Theorem~\ref{bigtheorem}. For the purpose of stating the theorem, let $\hat{f}_t(z)=f_t(z+V_t)$, where $f_t: \mathbb{H} \to H_t$ is the chordal $\SLE_\kappa$ Loewner chain and $V_t$ is the Loewner driving function for $f_t$. (See Section \ref{SLEsec} for definitions.) Further, let
\[
 \rho(\beta) = \frac{\kappa}{8(\beta+1)}
  \, \left[\left(\frac {\k + 4}\k\right) \, (\beta + 1) -
  1 \right]^2,
\]
and set
\[
\beta_{\pm}=-1 + \frac{\kappa}{12 + \kappa \mp
    4 \sqrt{8 + \kappa}}.
\]
Define 
\[
\Theta_\beta = \{t \in (0,2]:
     y|\hat f_t'(iy)| \wapprox y^{1-\beta} \}.
\]
(See Section~\ref{notsec} for the definition of $\wapprox$.)
\begin{theorem*}[Tip multifractal spectrum]
Suppose $\kappa>0$ and that $\beta_- \leq \beta \leq \beta_+$. For chordal $\SLE_\kappa$, with probability one,
\begin{equation*}
    \hdim(\Theta_\beta)  =   \frac{2 - \rho(\beta)}{2} , \quad \hdim[\gamma(\Theta_\beta) ]  
=  \frac{2-\rho(\beta)}{1-\beta},
\end{equation*}
\end{theorem*}
See the precise statement in Theorem~\ref{bigtheorem}. (We prove more than we state here.)

Notice that we obtain Beffara's theorem on the dimension of the $\SLE_\kappa$ path \cite{Beffara} as a corollary of Theorem~\ref{bigtheorem}.

Using the tip multifractal spectrum and some additional work we can derive the almost sure spectrum for harmonic measure at the tip; see Section~\ref{hmsec} for more details. 
Although we modify the definition of the spectrum somewhat, we prove in Theorem~\ref{harmtheorem} the stronger almost sure version of the theorem. To state it, define 
\[
\alpha_{\pm}=\frac{1}{1-\beta_{\pm}},\]
where $\beta_{\pm}$ are as above. 
\begin{theorem*}[Multifractal spectrum for harmonic measure at the tip] \label{hm_formula}
Suppose $\kappa>0$ and that $\alpha_- \le \alpha \le \alpha_+$. For chordal $\SLE_\kappa$, with probability one,
\begin{equation}  
\hdim[\gamma(\Theta_\alpha^\hm)]  = \alpha\left(1-\fr{4}{\kappa}\right) + \fr{(4+\kappa)^2}{8 \kappa}-\fr{\kappa}{8} \left( \frac{\alpha^2}{2\alpha -1} \right).
\end{equation}
\end{theorem*}  
In the final section of the paper we prove Theorem~\ref{newtheorem} which together with Theorem~\ref{harmtheorem} and a Beurling estimate shows that the right hand side of \eqref{hm_formula} gives the harmonic measure spectrum for a (one-sided) version which is closer to the usual definition, but for a smaller range of $\alpha$. 
 
\begin{figure}[t]
\centering
\includegraphics[width=85 mm]{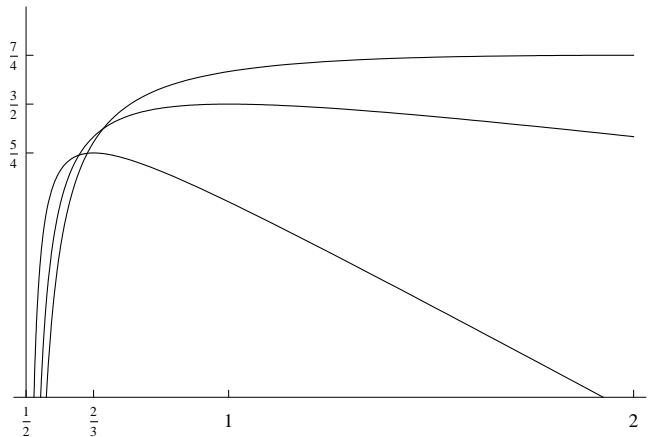}
\caption{Multifractal spectrum of harmonic measure at the tip for 
$\SLE_\kappa$, $\kappa =2,4,6$. The maximum is
the Hausdorff dimension of the curve.}
\end{figure}

\subsection{Outline of the paper}
Our paper is organized as follows. The next section discusses
some preliminary facts.  After setting some notation about asymptotics
in Section \ref{notsec}, the deterministic Loewner equation is discussed
in Section \ref{detersec}.  Much in this subsection
is standard but we have included
this in order to phrase the results appropriately for our purposes.
Also, we want to separate estimates that deal only with the Loewner
equation itself from those that are particular to $\SLE$.
In this subsection, there are three kinds of results: those that
hold for all conformal maps of $\Half$ for which we use the letter
$h$; those that hold for all solutions of the chordal Loewner equation
for which we use $g_t$ and $f_t = g_t^{-1}$; and towards the end fact
about solutions of the Loewner equation for driving functions that
are weakly H\"older-($1/2$). We also formally define the tip multifractal
spectra in this section. 

  The main theorem is not stated in full until
Section \ref{SLEsec} where the Schramm-Loewner evolution ($\SLE$),
that is, the solution of the Loewner equation with a Brownian motion input,
is discussed.  From here on a value of the $\SLE$ parameter $\kappa$ is
fixed and a large number of $\kappa$-dependent parameters are
defined.   Although we do not discuss it directly, what we are doing
is establishing the guess for the value of the multifractal spectrum
in terms of the Legendre transform of a logarithmic moment
generating function.

The basic proofs of the main theorem can be found in
Section \ref{proofsec}.  This section is relatively short because it
relies on estimates on the moments of the derivative some of
which were established in \cite{Law1,JL}; the necessary additions
are proved in Section \ref{estsec}.   Section \ref{newsec}
uses the {\em forward} Loewner flow to prove a result on the
harmonic measure spectrum stated in Section \ref{harmonsec}.  We
warn the reader that some of 
the notation in Section \ref{newsec} does not
agree with that earlier and that the assumption $\kappa < 8$
is made there.

\section{Preliminaries}

\subsection{Notations}   \label{notsec}

 In order to avoid writing bulky expressions with ratios
of logarithms, we will adopt the following notations.

We call a
 function $\psi:[0,\infty) \rightarrow (0,\infty)$ 
a {\em (positive)
subpower function} if it is continuous, nondecreasing,
and
\[          \lim_{x \rightarrow \infty} x^{-u}\, \psi(x) = 0, \]
for all $u > 0$.  


If $f$ and $g$ are positive functions tending to zero with $y$, we write
\begin{equation}  \label{nov18.1}
  f(y) \wapprox g(y) , \;\;\;\;  y \rightarrow 0+ , 
\end{equation}
if there exists a subpower function $\psi$ such that
\[      \psi(1/y)^{-1} \, g(y) \leq f(y) \leq \psi(1/y) \, g(y) ,
\quad y \to 0+ \]
We write 
\[    f(y) \lap g(y), \;\;\;\;  y \rightarrow 0+, \]
if
\[    \limsup_{y\rightarrow 0+} \frac{\log g(y)}
   {\log f(y)} \leq 1, \]
and 
\[        f(y) \lapp g(y), \;\;\;\;  y \rightarrow 0+, \]
if
\[    \liminf_{y\rightarrow 0+} \frac{\log g(y)}
   {\log f(y)} \leq 1. \]
   
Here i.o. stands for ``infinitely often''. Clearly
$f(y) \lap g(y)$ implies $f(y) \lapp g(y)$ but the converse
is not true. Similarly we write 
$f(y) \gap g(y)$ and $f(y) \gapp g(y)$ for
\[   \liminf_{y\rightarrow 0+} \frac{\log g(y)}
   {\log f(y)} \geq 1 \;\;\;\mbox{ and } \;\;\;
\limsup_{y\rightarrow 0+} \frac{\log g(y)}
   {\log f(y)} \geq 1 , \]
respectively.  
We write $f(y) \approx g(y)$ if  $f(y) \lap g(y)$ and
$f(y) \gap g(y)$, that is, if 
\[      \lim_{y \rightarrow 0+} \frac{\log g(y)}
   {\log f(y)} = 1. \]
Note that  
if $\beta > 0$, then
\[ f(y) \approx y^{\beta} \;\;\;\; \Longleftrightarrow  
\;\;\;\; f(y) \wapprox y^{\beta}. \]


 We  will also use the notations
 for
asymptotics for functions $f(n), g(n)$ as $n \rightarrow \infty$ along the positive integers.   

\subsection{Chordal Loewner equation}  \label{detersec}

In this section, we review some facts about conformal mappings
and the chordal Loewner equation.  See \cite[Chapters 3,4]{LConv}
for proofs of theorems stated without proof here. 

Suppose $\tilde \gamma: (-\infty, \infty) \to  \mathbb{C}$ is a curve as in the introduction.
The chordal Loewner equation is an equation that describes 
the evolution of $\tilde \gamma(0,\infty)$ given $\tilde \gamma(-\infty,0]$.
Let  $\tilde g$ be a conformal transformation
of $ \C \setminus \tilde \gamma(-\infty,0]$ onto the upper
half plane $\Half$ with $\tilde g(\gamma(0)) = 0$, 
$\tilde g(\infty) = \infty$.  In order to describe
$\tilde \gamma(t), t > 0$, it suffices to describe
\[    \gamma(t) := \tilde g(\tilde \gamma(t)) , \;\;\;\;
   0 \leq t < \infty , \]
and this is what the Loewner equation in $\Half$ does.
For the remainder of the paper, we will consider a curve $\gamma$
in $\Half$ as above.
The Riemann mapping theorem implies that there is
a unique conformal transformation
 $g_t$  of $\Half \setminus
\gamma(0,t]$ onto $\Half$ with $g_t(z) = z + o(1)$
as $z \rightarrow \infty$. We can expand
$g_t$
 at infinity, 
\[             g_t(z) = z + \frac{a(t)}{z} + O(|z|^{-2}) , \]
where $a(t)$ by definition is the half-plane capacity of $\gamma(0,t]$.
It is continuous and strictly increasing.  We make the 
(slightly) stronger assumption that $a(t) \rightarrow \infty$
as $t \rightarrow \infty$.   Then the chordal  
Loewner integral equation
states that
\[          g_t(z) = z + \int_0^t \frac{da(s)}
                        {g_s(z) - V_s} , \;\;\;\;\;
   t \leq T_z , \]
where $V_s = g_s(\gamma(s))$
and  $T_z = \inf\{t: \Im[g_t(z)] = 0 \} = 
\inf\{t: g_t(z) - V_t = 0\}$.  It can be shown that
$s \mapsto V_s$ is a continuous function.
  It is convenient to choose
a parameterization of $\gamma$ such that $a(t) = at$ for
some $a > 0$ in which case we get the Loewner
 differential equation
\begin{equation}  \label{chordal2}
      \p_tg_t(z) = \frac{a}{g_t(z) -V_t}, \;\;\;\;
    g_0(z) = z . 
\end{equation}
Let 
\[   f_t(z) = g_t^{-1}(z), \;\;\;\;
  \hat f_t(z) = f_t(z + V_t) = 
    g_t^{-1}(z + V_t) . \]
By differentiating both sides of
$f_t(g_t(z)) = z$ with respect to $t$ we see that ($f_t'(z):=\partial_z f_t(z)$)
\begin{equation}  \label{inversechord}
     \p_tf_t(z) = -  f_t'(z) \, \frac{a}{z - V_t},
\end{equation}
and since $g_t(\gamma(t)) = V_t$, we get
\begin{equation}  \label{curvedef}
    \gamma(t) = f_t(V_t) =
\lim_{y \rightarrow 0+} f_t(V_t + iy)
    = \lim_{y \rightarrow 0+} \hat f_t(iy). 
\end{equation}
We let
\[   v_t(y) = v(\hat f_t;y) = \int_0^y
   |\hat f_t'(iu)| \, du. \]
As mentioned before, if $v_t(y) < \infty$ for some $y > 0$,
then $v_t(0+) = 0$ and the limit in \eqref{curvedef} exists.
More work is needed to determine whether or
not $\gamma$ is
a continuous function of $t$. 
Note that if $g_t$ satisfies \eqref{chordal2} and $ g_t^*=
g_{t/a}$, then 
\[  \p_t g_t^*(z) = \frac{1}{g_t^*(z) - V_t^*}, \;\;\;\;
    g_0^* (z) = z, \]
where $V_t^* = V_{t/a}$.

Conversely, we can start with a continuous function
$t \mapsto V_t$ and $a > 0$ and define a Loewner chain $(g_t, \, t \ge 0)$
by \eqref{chordal2}.  We define $\gamma(t)$ by \eqref{curvedef}
provided that the limit exists.
 We say that the family of conformal maps $g_t$
{\em generates a curve} if $\gamma$  
 exists and is a continuous
function of $t$.  We do {\em not} assume that the
curve is simple. If  $H_t$ denotes the unbounded
component of $\Half \setminus \gamma(0,t]$, then $g_t$
is the unique conformal transformation of $H_t$ onto
$\Half$ satisfying
\[   g_t(z) = z + \frac{at}{z} + O(|z|^{-2}), \;\;\;\;
  z \rightarrow \infty. \]

\begin{lemma}  For
  every $t$ and every $y > 0$ with $v_t(y) < \infty$, 
\begin{equation}  \label{nov14.1}
             \frac{y \, |\hat f_t'(iy)|}
  {4} \leq   |\gamma(t) - \hat f_t(iy)|
    \leq v_t(y). 
\end{equation}
\end{lemma}

\begin{proof}
  The second estimate is immediate from the definition
of $v_t(y)$ and the first inequality follows from the Koebe
$(1/4)$-theorem applied to $\hat f_t$ on the open disk of radius $y$ about
$iy$.
\end{proof}

\begin{lemma}  \label{lemma34}
 If
$f_t$ satisfies \eqref{inversechord} and  $z = x+iy\in \Half$, then
for  $s \geq 0$
\begin{equation}  \label{nov19.9}
           e^{-5as/y^2}\, |f_t'(z)| 
 \leq  |f_{t+s}'(z)|\leq
       e^{5as/y^2}\, |f_t'(z)|  . 
\end{equation}
In particular, if $s \leq y^2$,
\[            e^{-5a}\, |f_t'(z)| 
 \leq  |f_{t+s}'(z)|\leq
       e^{5a}\, |f_t'(z)|  . \]
\end{lemma}

\begin{proof} Without loss of generality we may
assume that
$a=1$.  Differentiating \eqref{inversechord} with respect to $z$ yields
\[    \p_tf_t'(z)    = -f_t''(z) \, \frac{1}{z-V_t}
            + f_t'(z) \, \frac{1}{(z-V_t)^2} . \]
Note that $|z-V_t| \geq y$. 
Applying Bieberbach's theorem (the $n=2$ case of
the Bieberbach conjecture) to the disk
of radius $y$ about $z$, we can see that
\[  |f_t''(z)| \leq  4 \, y^{-1} \, |f_t'(z)|. \]
and hence
\[      |\p_t f_t'(z)| \leq 5  \, y^{-2} \, |f_t'(z)|,\]
which implies \eqref{nov19.9}.
\end{proof}

 The Koebe distortion and growth theorems are traditionally
stated in terms of univalent functions defined on
the unit disk, see, e.g., \cite[Chapter 2]{Pommerenke}. We will use these theorems for univalent
functions on $\Half$, and the next proposition gives
the appropriate results.

\begin{proposition} \label{prop34}
 Suppose $h: \Half \rightarrow \C$ is a
conformal transformation, $ x \in \R, y > 0,$ and $r \geq 1$.
Then    
\begin{equation}  \label{nov19.10}
 (x^2 + 4)^{-2}  
 \, |h'(iy)| \leq  |h'(y(x+i))|  \leq     \, (x^2 + 4)^2
  \, |h'(iy)|,
\end{equation}
 \begin{equation}  \label{nov19.12}
  |h(y(x+i)) - h(iy) | \leq  \frac{(x^2+4)^{3/2} \, |x| }{2} \,
  y\, |h'(iy)|, 
\end{equation}
\begin{equation}  \label{nov19.11}
   r^{-3} \, |h'(iy)|  \leq
   |h'(iyr)| \leq r \, |h'(iy)| , 
\end{equation}
\begin{equation}  \label{nov19.13}
   | h(iyr) - h(iy)|  \leq \frac{r^2 - 1}{2}
 \, y \, |h'(iy)|.
\end{equation}
\end{proposition}

\begin{proof} By scaling  we may assume that $y=1$.
Let 
\[   G(z) = \frac{z-i}{z+i} , \;\;\;\;\;  G'(z) = \frac{2i}{(z+i)^2}, \]
which is a conformal transformation of $\Half$ onto
the unit disk $\Disk$  with
$G(i) = 0, |G'(i)| = 1/2$. 
 We can write
\[           h(z) =  f(G(z)) ,\;\;\;\;\;  h'(z) = f'(G(z))
  \, G'(z),  \]
where $f$ is a univalent function on $\Disk$. 
 The distortion theorem
tells us that
\[           {|f'(w)|} 
     \leq \frac{1+|w|}{(1-|w|)^3} \, |f'(0)| ,\;\;\;\;
  |w| < 1,
\]
and the growth theorem
states that
\[   |f(w) - f(0)| \leq \frac{|w|}{(1-|w|)^2}\, |f'(0)|, \;\;\;\;
  |w| < 1. \]
Since $|G'(i)|=1/2$,  we get
\begin{equation}  \label{nov19.111}   
      { |h'(z)|}\leq \frac{2\,|G'(z)|\,
 [1 +|G(z)|]}{ (1-|G(z)|)^3}\; {|h'(i)|} , 
\end{equation}
and 
\begin{equation}  \label{nov19.112}
|h(z) - h(i)| \leq  \frac{2\,|G(z)|}
                     {(1 - |G(z)|)^2} \, |h'(i)|. 
\end{equation}

Since
\[  |G(x+i)| = \frac{|x|}{\sqrt{x^2+4}} , \;\;\;\;\;
    |G'(x+i)| = \frac{2}{x^2+4}, \]
we plug into \eqref{nov19.111} and see that 
\[    |h'(x+1)| \leq \frac 1{16} \,
[{\sqrt{x^2 + 4} + |x|}]^4 \, |h'(i)|
  \leq   (x^2+4)^2 \, |h'(i)|. \]
This gives the second inequality in \eqref{nov19.10} and
the first follows easily by real translation.
Plugging into \eqref{nov19.112} gives
\begin{eqnarray*}
 |h(x+i) - h(x)| & \leq & \frac{2\, |x| \, \sqrt{x^2 + 4}
  \, (\sqrt{x^2 + 4} + x)^2}{16} \, |h'(i)| \\
 & \leq &
              \frac{(x^2 + 4)^{3/2} \, |x|}{2} \, |h'(i)|. 
\end{eqnarray*}

Since  $r \geq 1$,
\[    |G(ir)| = \frac{r-1}{1+r} = |G(i/r)|, \;\;\;\;\;
    |G'(ir)| = \frac{2}{(1+r)^2}, \;\;\;\;
  |G'(i/r)| = \frac{2r^2}{(r+1)^2}. \]
Plugging into \eqref{nov19.111}  and \eqref{nov19.112}
gives  \eqref{nov19.11} and \eqref{nov19.13}.

\end{proof}

\begin{corollary}  \label{apr10.cor1}
If $h: \Half \rightarrow \C$ is a conformal transformation,
then for every $y >0$,
\begin{equation} \label{vest}
      v(h;y) \geq \frac {y\, |h'(iy)|}2, 
\end{equation}
\[         \frac 23\,   v(h;2^{-n})
  \leq             \sum_{j=n}^\infty
                 2^{-j}\,  |h'(i2^{-j})|    \leq   
\frac {8}3 \, v(h;2^{-n}). \]
\end{corollary}

\begin{proof}  
We write
\[   v(h;2^{-n}) = \sum_{j=n}^\infty
    \int_{2^{-j-1}}^{2^{-j}}  |h'(iy)| \, dy . \]
Using \eqref{nov19.11} (which holds for $r > 1$)
\[ v(h;y) = \int_0^y |h'(is)| \, ds =
   y \int_0^1 |h'(iry)| \, dr \geq y \, |h'(iy)|
  \, \int_0^1 r \, dr =  \frac {y\, |h'(iy)|}2,  \]
\[ 
\int_{r/2}^{r} |h'(iy)| \, dy  
  \leq r \, 
|h_t'(ir)| \int_{1/2}^1 s^{-3} \, ds 
  =  \frac{3r}{2} \, |h'(ir)| ,\]
\[ \int_{r/2}^{r} |h'(iy)| \, dy  
    \geq r \, |h'(ir)| \int_{1/2}^1 s \, ds
  = \frac{3r}{8} \, |h'(ir)|. \]
\end{proof}

 We define the following measure of the
modulus of continuity of $U_t$:
\[  \Delta(t,s) = \sup_{0 \leq r \leq s^2}
    \sqrt{s^{-2} (U_{t+r} - U_t)^2 + 4}. \]
Note that $\Delta(t,s) \geq 2$ and it  is
of order one if 
\[  \sup_{0 \leq r \leq s^2} |U_{t+r} - U_t| \approx s.\] The
 definition of $\Delta(t,s)$
with the $4$ has been chosen to make
the statement of the next proposition cleaner.

\begin{proposition}  \label{prop.456}
If $t \geq 0$ and $0 \leq r \leq s^2$ with
$v_t(s) + v_{t+r}(s) < \infty$, then
\begin{equation}  \label{nov20.2}
 |\gamma(t+r) - \gamma(t)|
   \leq v_t(s) + v_{t+r}(s) + 
e^{5a} \, |\hat f_t'(si)|
  \, \Delta(t,s)^4\,s. 
\end{equation}
\end{proposition}

\begin{proof} By the triangle inequality and
\eqref{nov14.1},
\begin{eqnarray*}
\lefteqn{ |\gamma(t+r) - \gamma(t)| } \hspace{.3in} \\
& \leq & |\hat f_t(si) - \gamma(t)|+
  |\hat f_{t+r}(si) - \gamma(t+r)| + |\hat f_{t}(si) - \hat f_{t+r}(si)| 
  \\
&  \leq  & v_t(s) + v_{t+r}(s) + 
  |f_{t+r}(V_{t+r} + si) - f_t(V_t + si)|.
\end{eqnarray*}
Also,
\[  |f_{t+r}(V_{t+r} + si) - f_t(V_t + si)| \leq \hspace{2.5in} \]
\[ 
   |f_{t+r}(V_{t+r} + si) - f_{t+r}(V_t + si)|
    + |f_{t+r}(V_t + si)- f_t(V_t + si)|.\]
Using \eqref{nov19.12} and \eqref{nov19.9},  
we see that 
\begin{eqnarray*}   |f_{t+r}(V_{t+r} + si) - f_{t+r}(V_t + si)|
  & \leq & \frac 12 \, |f_{t+r}'(V_t + si)| \, \Delta(t,s)^4\,s \\
  & \leq &  \frac 12 \,
e^{5a} \, |f_t'(V_t+si)|\, \Delta(t,s)^4\,s\\
  &= & \frac 12 \,
e^{5a} \, |\hat f_t'(si)|\, \Delta(t,s)^4\,s
\end{eqnarray*}
Also \eqref{nov19.9} and \eqref{inversechord} imply that
\begin{equation}  \label{nov20.19}
 |\p_r f_{t+r}(V_t+ is)| \leq \frac as \, e^{5ar/s^2} \,
  |\hat{f}_t'(is)|
 , 
\end{equation}
and hence
\begin{eqnarray}   |f_{t+r}(V_t+ is) - f_t(V_t +is)| & \leq& 
            |f_t'(V_t+ is)|  \int_0^{s^2} \frac as \,
 e^{5au/s^2} \, du  \nonumber\\
 & = &  \frac{s}{5} \, e^{5a} \, |\hat f_t'(is)|\\
 &  < & \frac{s}{2} \, e^{5a} \, \Delta(t,s)^4 \, |\hat f_t'(is)|. 
\label{nov20.20}
\end{eqnarray}

\end{proof}

\begin{lemma}  \label{ber}
There exist $c > 0$ such that for $t \geq 0$ and $0 < y \leq 1$,
\[    c \,\frac{y}{\sqrt{2at+1}}
\leq |\hat f_t'(iy)| \leq   \frac{\sqrt{2at + 1}}
  y  . \]
\end{lemma}

\begin{proof}  We may assume $a=1$ for otherwise we consider $ g_t^*
= g_{t/a}$.   Let $w= \hat f_t(iy)$, that is, $g_t(w) = V_t + iy$, and
let $Y_t = \Im[g_t(w)]$.  The
Loewner equation implies that $\p_t [Y_t^2 ]\geq -2$ and hence
$\Im[w] \leq  \sqrt{2t+ 1}$.  Similarly, $\Im[\gamma(s)] \leq 
\sqrt{2t} \leq \sqrt{2t + 1}$ for $0 \leq s \leq t $.  The Loewner
equation also implies that $\p_t[Y_t/|g_t'(w)|] \leq 0$, which
implies
\[  y \, |\hat f_t'(iy)| = \frac{Y_t}{|g_t'(w)|}
   \leq \frac{Y_0}{|g_0'(w)|} = \Im[w] \leq \sqrt{2t + 1}. \]
This gives the second inequality. 

For the first inequality, let
$d = \dist(w,\gamma[0,t] \cup \R)$.  The Beurling estimate
\cite[Theorem 3.76]{LConv}
implies that there is a $c_* < \infty$ such that the probability
that a Brownian motion starting at $w$ goes distance $\sqrt{2t+1}$
without hitting $\gamma[0,t] \cup \R$ is bounded above by
\[               c_* \, \left(\frac{d}{\sqrt{2t+1}}\right)^{1/2}. \]
By the gambler's ruin estimate, the probability that a Brownian motion
in $\Half$
starting at $iy$ reaches $I_t :=\{\tilde w: \Im[\tilde w]
  = 2 \sqrt{2t +1}\}$ before hitting the real line equals
$y \,[2\sqrt{2t+1}]^{-1} . $  Since the imaginary part decreases
in the forward Loewner flow, it follows from conformal invariance that the probability
that a Brownian motion starting at $w$ reaches $I_t$ before hitting
$\gamma[0,t] \cup \R$ is at least $y \,[2\sqrt{2t+1}]^{-1} . $  Therefore
\[                 c_* \, \left(\frac{d}{\sqrt{2t+1}}\right)^{1/2}
     \geq \frac{y}{2\sqrt{2t+1}}. \]
The Koebe $(1/4)$-theorem implies that $d \leq 4 \, y \, |\hat f_t'(iy)|$,
and plugging in we get
\[                  |\hat f_t'(iy)| \geq \frac{y}{16 \, c_*^2
  \, \sqrt{2t+1}}. \]
\end{proof}

\begin{proposition}\label{standard2}
Suppose $h: \Half \rightarrow \C$
is a conformal transformation and
$v(h;y)$ is as defined in \eqref{nov18.12}.
  Then for every $\beta < 1$, as $y \rightarrow 0+$
\[ y|h'(iy)| \lap y^{1-\beta} \mbox{ if and only if }
    v(h;y) \lap y^{1-\beta}, \]
\[ y|h'(iy)| \gapp y^{1-\beta} \mbox{ if and only if }
          v(h;y) \gapp y^{1-\beta}. \]
\begin{equation}  \label{nov18.7}
        y|h'(iy)| \wapprox y^{1-\beta }
\mbox{ if and only if }
       v(h;y)  
 \wapprox y^{1-\beta}.
\end{equation}
\end{proposition}

\begin{proof} Using Corollary \ref{apr10.cor1}, all of
the assertions follow easily except  that
$v(h;y) \wapprox y^{1-\beta}$ implies $y|h'(iy)|
 \wapprox y^{1-\beta}$ which we will show here.
Assume $v(h;y) \wapprox y^{1-\beta}$.
  By \eqref{vest}, we know that
$y|h'(iy)| \lap y^{1-\beta}.$ For  $0 < \epsilon < 1-\beta$, 
let
\[       
   u = u_\epsilon = {\frac{3\epsilon}{1-\beta-\epsilon}}, \]
and note that \eqref{nov19.11} implies for $y$ sufficiently
small
\begin{eqnarray*}
      y^{1-\beta + \epsilon}    \leq  
  v(h;y)  
  & =  &  v(h;y^{1+u}) + \int_{y^{1+u}}^y |h'(is)| \, ds \\
& 
  \leq   & y^{1-\beta + 2 \epsilon} +  \int_{ y^{1+u}}^y |h'(is)| \, ds 
 \\ &  \leq & y^{1-\beta + 2 \epsilon} +   
  y^{1-3u
}\, |h'(iy)|.
\end{eqnarray*}
Hence for all $y$ sufficiently small,
\[   y|h'(iy)| \geq \frac{1}{2} \, y^{1-\beta} \, y^{3u_\epsilon + 
\epsilon}. \]
Since $u_\epsilon \rightarrow 0$ as $\epsilon \rightarrow 0+$, this
gives
$           y|h'(iy)| \gap y^{1-\beta}. $
\end{proof}

\begin{definition}
For every $-1 \leq \beta \le 1$, let
\[  \overline \Theta_\beta = \{t \in (0,2]:
    y|\hat f_t'(iy)| \gapp y^{1-\beta}\},  \]
\[  \Theta_\beta = \{t \in (0,2]:
     y|\hat f_t'(iy)| \wapprox y^{1-\beta} \}, \]
\[  \tilde \Theta_\beta = \{t \in (0,2]:
      y|\hat f_t'(iy)| \lap y^{1-\beta} \} , \]
\[  \underline \Theta_\beta = \{t \in (0,2]:
      y|\hat f_t'(iy)| \lapp y^{1-\beta} \} , \]
\[    \underline \Theta_\beta^* =
    \{t \in (0,2]:
      v_t(y) \lapp y^{1-\beta} \} , \]
where in each case the asymptotics are as
$y \rightarrow 0+$.
\end{definition}

If $\beta \neq 1$, we can write these sets as the
set of $t \in (0,2]$ such that
\[    \limsup_{y \rightarrow 0+} \frac{\log
   |\hat f_t'(iy)|}{\log (1/y)} \geq \beta, \]
\[   \lim_{y \rightarrow 0+} \frac{\log
   |\hat f_t'(iy)|}{\log (1/y)}  = \beta, \]
\[    \limsup_{y \rightarrow 0+} \frac{\log
   |\hat f_t'(iy)|}{\log (1/y)} \leq \beta, \]
\[    \liminf_{y \rightarrow 0+} \frac{\log
   |\hat f_t'(iy)|}{\log (1/y)} \leq \beta, \]
 \[    \liminf_{y \rightarrow 0+} \frac{\log
   v_t(y)}{\log (1/y)} \leq \beta-1, \]
respectively.

Using Lemma \ref{ber}, we can
see that 
  for every $\beta > 1$,
\[   \overline \Theta_\beta = \Theta_\beta
   = \Theta_{-\beta} = \tilde \Theta_{-\beta} =
 \underline \Theta_{-\beta} = \emptyset.\]
 Note that  \eqref{vest} implies that $\underline \Theta_\beta^*
\subset \underline \Theta_\beta$.
Using Proposition \ref{standard2} we can see that
we can also write
\[  \Theta_\beta = \{t\in(0,2]: v_t(y) \wapprox
    y^{1-\beta}, \;\; y \rightarrow 0+\}, \]
and similarly for $\overline \Theta_\beta, \tilde
\Theta_\beta$.
 Also, $\overline \Theta_\beta \cup \tilde \Theta_\beta
= (0,2]$ and
\[               \Theta_\beta \subset \overline \Theta_\beta
   \cap \tilde \Theta_\beta \cap \underline \Theta_\beta^*.\]

\begin{definition}
The driving function $V_t$ is {\em weakly H\"older-(1/2)}
 on $[0,2]$ if for  
each $\alpha < 1/2$, $V_t$ is H\"older continuous
of order $\alpha$ on [0,2].
\end{definition}  

Two equivalent definitions are
\begin{itemize}
\item If
\[ \delta(s) = \sup \{|V_{t+s} - V_t|:0 \leq t \leq t+s \leq 2
  \}, \]
then
\[ \psi(x) = \sup_{s \geq 1/x }\;  
 s^{-1/2} \, \delta(s) \]
 is a subpower function.
\item 
There is a subpower function
$\psi$ such that for all $0 \leq t \leq 2, 0 \leq s \leq 1$,
\[            \Delta(t,s) \leq \psi(1/s). \]
\end{itemize}
The next proposition shows that for
weakly H\"older-($1/2$) functions $V_t$, it suffices 
to consider    
dyadic times  in the definition of $\Theta_\beta$, etc.

\begin{proposition}   \label{distortprop2}
 Suppose $V_t$ is weakly H\"older-($1/2$) on
$[0,2].$  For each $t \in [0,2]$ define $t_n
= t_n(t) 
= (j-1)/2^{2n}$ if
\[            \frac{j-1}{2^{2n}} \leq t <
    \frac{j}{2^{2n}}. \]
Then for $-1 \leq \beta \le 1$, 
\begin{itemize}
\item   
\[  \Theta_\beta   = 
\left\{t \in (0,2]:   
 { 2^{-n}|\hat f_{t_n}'(i2^{-n})|}  \wapprox
   2^{-n(1-\beta)} \right\},
\]
\[ \overline
 \Theta_{\beta} =
\left\{t\in (0,2]:   
    2^{-n}|\hat{f}_{t_n}'(i2^{-n})|   \gapp 2^{-n(1-\beta)}  \right\},
\]
\[ \underline
 \Theta_{\beta}  =
\left\{t\in (0,2]:   
    2^{-n}|\hat{f}_{t_n}'(i2^{-n})|  \lap_{{\rm i.o}} 2^{-n(1-\beta)}  \right\},
\]
where the asymptotics are as $n \rightarrow \infty$ along
the integers. 

\item  If $t \in 
 \overline \Theta_{\beta}$,
\[   { v_t(y)}  \gapp y^{1-\beta}, \;\;\;\;
  |\gamma(t) - \hat f_t(iy)| \gapp
          y^{1-\beta}   ,\;\;\;\;  y \rightarrow 0+.\]

\item If $t \in \tilde  \Theta_{\beta}$,
\begin{equation} \label{110}
  { v_t(y)}  \lap y^{1-\beta}, \;\;\;\;
  |\gamma(t) - \hat f_t(iy)| \lap
          y^{1-\beta}   ,\;\;\;\;  y \rightarrow 0+.
\end{equation}

\item  If $t \in \Theta_\beta$, 
\[   v_t(y) \wapprox y^{1-\beta} , \;\;\;\;
   |\gamma(t) - \hat f_t(iy)| \wapprox y^{1-\beta},
\;\;\;\; y \rightarrow 0+. \]

\end{itemize}
\end{proposition}

\begin{proof}
Note that 
\begin{eqnarray*}
 |\hat f_{t}'(i2^{-n})| = |f_{t}'(V_{t} + i2^{-n})|
    & \leq  & \Delta(t,2^{-n})^4 \, |f_{t}'(
    V_{t_n} + i2^{-n})|\\
&  \leq & e^{5a} \, 
  \Delta(t,2^{-n})^4 \, |\hat f_{t_n}'(i2^{-n})|, 
\end{eqnarray*}
and similarly
\[       |\hat f_{t}'(i2^{-n})| \geq e^{-5a} \,
 \Delta(t,2^{-n})^{-4} \, |\hat f_{t_n}'(i2^{-n})|. \]
 Hence if $V_t$ is weakly H\"older$-(1/2),$
 there is a subpower function
$\psi$ such that for all $t,n$,
\[   \psi(2^{n})^{-1} \, |\hat f_{t_n}'(i2^{-n})|
          \leq |\hat f_t'(i2^{-n})|
  \leq   \psi(2^{n})\, |\hat f_{t_n}'(i2^{-n})|.\]
This implies the first assertion and the remaining follow from 
\eqref{nov14.1}. 
  
\end{proof}

\subsection{Harmonic measure at the tip}\label{hmsec}

We will now discuss harmonic measure giving two 
nonequivalent definitions,
one that is standard and one which is more directly related to
the multifractal spectrum we have discussed.

In this subsection $\gamma$ denotes a curve in $\HH$ with
one endpoint on the real line. We assume that the curve comes from a
Loewner chain driven by a continuous function
$V_t$, so it may have double points but it does not cross itself. 
Let $H_t$ be the unbounded connected component of $\HH \setminus
\gamma[0,t]$. As above, we write $g_t:H_t \to \HH$ for the normalized
conformal mapping so that $\lim_{y \to 0+}\hat{f}_t(iy)=\gamma(t)$,
where $f_t=g_t^{-1}$ and $\hat f_t(z)=f_t(z+V_t)$. If
the curve has double points, we are interpeting
$\gamma(t)$  in terms of prime ends, and
we then tacitly understand $\gamma(t)$ as the prime end corresponding to $V_t$.

If $z \in H_t$, then $\hm_{t,z}$ will denote the usual harmonic
measure of $\R \cup \gamma(0,t]$ starting at $z$, that is to say
the hitting measure of Brownian motion starting at $z$ stopped
when it reaches $\p H_t$.  We let
\[          
\hm_t(U) = \lim_{y \rightarrow \infty}
             y \, \hm_{t,iy}(U), 
\]
which is the normalized harmonic measure from the boundary
point at infinity.
Note that for each $z \in H_t$, $\hm_t$ and $\hm_{t,z}$ are
mutually absolutely continuous.  Also, conformal invariance, the normalization at infinity, 
and the well-known
Poisson kernel in $\Half$ show that for bounded
$U$, 
\[         \hm_t(U) = \frac{1}{\pi} \, {\rm length}
          \left[g_t(U)\right]. \]
Let
\[   \tilde \mu(t,\epsilon) = \hm_t\left[\overline \ball(\gamma(t),\epsilon)
  \right], \]
where $\ball(z,\epsilon)$ denotes the open disk of radius 
$\epsilon$ about $z$ with closure $\overline \ball(z,\epsilon)$.
For $\alpha > 0$, define
\[ \tilde 
\Theta_{\alpha}^{\hm} = \{t \in (0,2]: \tilde \mu(t,\epsilon) 
\wapprox \epsilon^{\alpha}, \, \epsilon \to 0+\}.\]
We define
the multifractal spectrum of harmonic measure at the tip by 
\[  \alpha \longmapsto \hdim\left[ \gamma(\tilde 
\Theta_{\alpha}^{\hm}) \right]. \]
This multifractal spectrum can be hard to compute.  One of the
difficulties is that $\ball(\gamma(t),\epsilon) \cap H_t$
can contain many connected components whose images under $g_t$
are far apart.  We  will give a different definition that
is more directly related to the tip multifractal spectrum
in this paper.

\begin{figure}[t]
\centering
\includegraphics[width=85 mm]{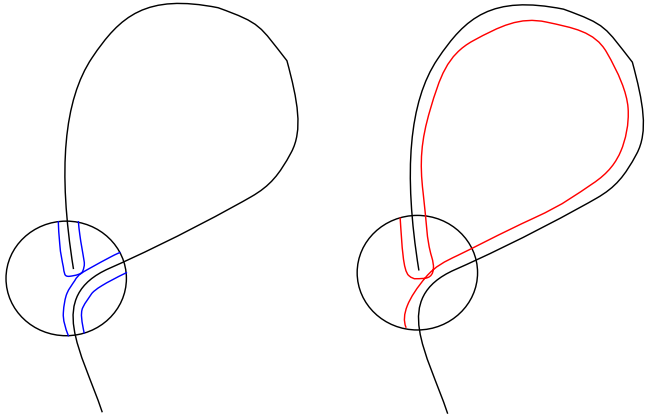}
\caption{Sketch showing the two different subsets of the boundary that we use to define $\tilde{\Theta}^{\rm{hm}}_{\alpha}$ (left) and $\Theta^{\rm{hm}}_{\alpha}$ (right).}
\end{figure}

Fix $t > 0, \epsilon > 0$ and let $\ball =
\ball(\gamma(t),\epsilon)$.   Let $O = O _{t,\epsilon}
$ denote the connected
component of $\ball \cap H_t $ that contains $\gamma(t)$
(considered as a prime end) on its boundary.  There
is a connected component (open arc), $\sigma =
\sigma_{t,\epsilon}$ 
 of $\p \ball \cap H_t$ that is in $\p O$ and such
that every curve from $\gamma(t)$ (again viewed as a prime
end) to infinity in $H_t$ passes through $\sigma$. We require also that $\sigma$ separates all other such open arcs from infinity in $H_t$; $\sigma$ is then unique.
 
Let $x_- = x_{-,t,\epsilon} < V_t < x_+ = x_{+,t,\epsilon}$
denote the images of the endpoints of $\sigma$ under $g_t$ (which always exist, see, e.g., \cite{Pommerenke}). 
Let $E = E_{t,\epsilon}$ denote the 
preimage of $[x_-,x_+]$ under $g_t$, and let
\[  \mu(t,\epsilon) = \hm_t[E] = \frac{x_+ - x_-}{\pi}. \]
It is not necessarily 
true that $E  \subset \overline \ball(\gamma(t),\epsilon)$.  
However an estimate using the Beurling projection theorem
shows that there is a $c< \infty$ such that 
\begin{equation}   \label{best2}
             \mu(t,\epsilon/2) \leq c \, 
\tilde \mu(t,\epsilon).
\end{equation}
   We define
\[
\Theta_{\alpha}^{\hm}=\{t \in (0,2]: \mu(t,\epsilon)
\wapprox \epsilon^{\alpha}, \, \epsilon \to 0+\}.
\]
The next lemma makes the connection with the tip multifractal spectrum.
\begin{lemma}  \label{aug13.lemma}
If  $1/2 \le \alpha < \iy$, then
\[
\Theta^{\hm}_{\alpha}=\Theta_{1-1/\alpha}.
\]
\end{lemma}
\begin{proof}
We will prove that there exist $0<c_1,c_2 < \iy$ such that for all $t \ge 0$ and all $\ee>0$ sufficiently small,
\begin{equation}  \label{aug13.1}
\mu(t,2v_t(\ee)) \geq  c_1 \ee,  
\end{equation}
\begin{equation} \label{aug13.2}
	\mu(t,\ee)|\hat{f}'_t(i\mu(t,\ee))| \le c_2 \ee.
\end{equation}
The lemma follows immediately from these estimates combined with Proposition~\ref{ber}.

Let $\eta_\ee$ denote the line segment
$(0,i\ee]$. The harmonic measure
from infinity of $\eta_\ee$
 in $\Half \setminus \eta_\ee$ equals
 $c_1\ee$ for a specific constant $c_1$, and hence by conformal invariance
the harmonic measure from infinity of $\eta_\ee^*:=
\hat f_t
\circ \eta_\ee$ in $H_t \setminus \eta_\ee^*$ is also
$c_1\ee$.  Since $\eta^*_\ee$ is a curve of length
$v_t(\ee)$ and one of its endpoints is $\gamma(t)$,
the interior of $\eta^*_\ee$ is contained in
$O_{t,v_t(\ee)}$. From this and a Beurling
estimate as in \eqref{best2}, 
we get \eqref{aug13.1}. 

It remains to prove \eqref{aug13.2}. To this end, let $\sigma_\epsilon = \sigma_{\epsilon,t}$ be the open arc whose endpoints are mapped
to $x_{-,\epsilon} < x_{+,\epsilon}$ as above. Let $\ell_{\ee}=x_{+,\epsilon}-x_{-,\epsilon}$ and note that $\mu(t, \ee)=\ell_{\ee}/\pi$. As $\epsilon \rightarrow 0+$, $\ell_{\ee} \rightarrow 0$ (using, e.g., the Beurling estimate), and hence for $y$ sufficiently small we can choose $\ee$ such that
$\ell_{\ee} \leq y$ and $\ell_{2\ee} \geq y$.  
Hence 
 it suffices to show that $2\epsilon  \geq c\, y \, |\hat{f}_t'(iy)|$.
Since $\ell_{2 \epsilon} \geq y$, we can
see that there exists $c_2$ such that  the probability
that a Brownian motion starting at $iy$ hits $g_t \circ
 \sigma_{2\epsilon}$ before leaving $\Half$ is at least
 $c_2$.  By conformal invariance, this is also true for
 a Brownian motion starting at $\hat f_t(iy)$ hitting
 $\sigma_{2\epsilon}$ before leaving $H_t$.  The
 distortion theorem and the 
 Koebe-(1/4) theorem show that $\dist[\hat f_t(iy),\p H_t]
 \asymp  \, y \, |\hat{f}_t'(iy)|$.
  The needed estimate then comes
 from the Beurling estimate which implies in any simply
 connected domain $D$, if $V \subset \p D$,
 \[   \hm_D(z,V) \leq c \, \sqrt{\frac{\diam(V)}{
              \dist(z,\p D)}} . \]
\end{proof}

\section{Tip spectrum for $\SLE$}  \label{SLEsec}

Let $\kappa > 0$ and $a = 2/\kappa$.   Then
the chordal Schramm-Loewner evolution with parameter $\kappa$
($\SLE_\kappa$)
is the solution to \eqref{chordal2} with $a = 2/\kappa$
where  $V_t$ is a standard Brownian motion.  
It is well known 
  that with probability one, $V_t$
is weakly H\"older-$(1/2)$. 
Let
\[   d = \min\left\{1 + \frac \k 8, 2 \right\}. \]
It was proved by Beffara \cite{Beffara} that $d$
is  the Hausdorff dimension of the path $\gamma[0,2]$.
This will follow as a particular case of our
main theorem,  so we will not need to
assume this result.  However,
it is convenient to use this notation.

 \subsection{Main theorem}

  Before stating the 
main theorem,  we will define some special values
of the parameter $\beta$.  See Section \ref{parasec}
for more details.
Let
\begin{equation}  \label{nov10.1}
 \rho(\beta) = \frac{\kappa}{8(\beta+1)}
  \, \left[\left(\frac {\k + 4}\k\right) \, (\beta + 1) -
  1 \right]^2, 
\end{equation}
\[  \hat d_\beta = \frac{2 - \rho(\beta)}{2},\]\[
   d_\beta =    \frac{2\, \hat d_\beta}{1-\beta}  
     = \frac{2-\rho(\beta)}{1-\beta}. \]
The maximum value of $\hat d_\beta$ equals $1$ and is
obtained at
\[      \beta_\# :=  \frac{\kappa}{\k + 4} - 1 . \]
The maximum value of $d_\beta$ equals $d$ and is
obtained at
\[ \beta_*   :=  \frac{\kappa}{\max\{4,\kappa - 4\}} - 1.\]
We define $\beta_- \leq \beta_\# \leq \beta_* \leq \beta_+$ by
$\rho(\beta_-) = \rho(\beta_+) = 2$.  A straightforward
computation gives
\begin{equation} \label{nov10.2}
  \beta_+   = -1 + \frac
  \k  {12 + \k - 4 \sqrt{8 + \k}}  , 
\end{equation}
\begin{equation} \label{nov10.3}
  \beta_- = -1 + \frac{\kappa}{12 + \kappa +
    4 \sqrt{8 + \kappa}} < 0.
\end{equation}
Also $-1 < \beta_- <\beta_+ \leq 1$ with equality
only for $\kappa = 8$.
\begin{remark}
The function $\beta_+(\kappa)$ determines the optimal H\"older exponent for the $\SLE_\kappa$ path in the capacity parameterization: With probability one, the chordal $\SLE_\kappa$ path away from the base is H\"older-$\alpha$ for $\alpha < (1-\beta_+)/2$ and not H\"older-$\alpha$ for $\alpha > (1-\beta_+)/2$. See Theorem~1.1 of \cite{JL} for a precise statement.
\end{remark}

\begin{theorem} \label{bigtheorem}
For chordal $\SLE_\kappa$, if $-1 \leq \beta \leq 1$, 
  the following
holds with probability one.

\begin{itemize}

\item   If $\beta_- \leq \beta \leq \beta_+,$ 
\begin{equation}  \label{sept18.0}
    \hdim(\Theta_\beta)  =   \hat  d_\beta , \;\;\;\;\;
 \hdim[\gamma(\Theta_\beta) ]  
=   d_\beta . 
\end{equation}

\item  If $\beta_\# \leq \beta \leq \beta_+,$  
\begin{equation}  \label{sept18.1alt}
  \hdim(\overline \Theta_\beta)  = \hat d_\beta .
\end{equation}

\item  If $\beta_* \leq \beta \leq \beta_+,$  
\begin{equation}  \label{sept18.1}
\hdim[\gamma(\overline \Theta_\beta)]  
   = d_\beta. 
\end{equation}

\item  
If $\beta_- \leq \beta \leq \beta_\#$, 
  \begin{equation}  \label{sept18.2alt}
 \hdim(\underline \Theta_\beta)  = \hat d_\beta .
\end{equation}

\item  
If $\beta_- \leq \beta \leq \beta_*$, 
  \begin{equation}  \label{sept18.2}
 \hdim[\gamma( \underline\Theta_\beta^*) ]  
   = d_\beta.  
\end{equation}

\item If $\beta >  \beta_+ $,
 $\overline \Theta_\beta   = \emptyset$.

\item  If $\beta < \beta_-$, $\underline \Theta_\beta
 = \emptyset$.

\end{itemize}

\end{theorem}

\subsection{Remarks}

\begin{itemize}

\item  It follows from the theorem that with probability one
the results hold for a dense set of $\beta$.  This implies
that with probability one, \eqref{sept18.1alt}--\eqref{sept18.2}
hold for all $\beta$. However, we have not
shown whether or not for a particular realization, there might
be an exceptional $\beta$ for which \eqref{sept18.0} does not hold.

\item The restriction to $t\in(0,2]$ is only a convenience.
By scaling we get a similar result for
$t \in (0,\infty)$.

\item  The relationship $\hdim[\gamma(\Theta_\beta)]
  = 2 \,\hdim[\Theta_\beta]/(1-\beta)$ can be
understood as follows.  For $s$ small, the image of the
interval $[t,t+s^2]$  under $\hat f_t$
can be approximated by a set of
diameter $s |\hat f_t'(is)|$ containing $\hat f_t
(is).$  If $|\hat f_t'(i s)| \approx s^{-\beta}$, then this set has
diameter $s^{1-\beta}$.  That is to say, intervals of length
(diameter) $s^2$ in a covering of $\Theta_\beta$
are sent to sets of diameter $s^{1-\beta}$.  Note that this is in
contrast to complex Brownian motion where intervals of length $s^2$
are always sent to sets whose diameter is of order $s$.

\item  Since $\Theta_\beta \subset \overline \Theta_\beta
\cap 
\tilde \Theta_\beta \cap  \underline \Theta_\beta^* $ 
and $\underline \Theta_\beta^* \subset \underline \Theta_\beta$,
it suffices to prove the lower bounds for
$\Theta_\beta$ in \eqref{sept18.0} 
and the upper bounds for $\overline
\Theta_\beta,  \underline \Theta_\beta, \underline
\Theta_\beta^*$ 
in \eqref{sept18.1alt}--\eqref{sept18.2}.
The upper bounds will be proved in Section \ref{uppersec}
and the lower bounds in Section \ref{lowersec}.

\item  To prove the upper bound \eqref{sept18.1alt} it suffices
to show for each $s > 0$,
\[           \hdim(\overline \Theta_\beta \cap (s,2])
  \leq \hat d_\beta, \]
and similarly for \eqref{sept18.1}--\eqref{sept18.2}.  This is what
we do in Section \ref{uppersec}.

\item  Recall that $\underline \Theta_\beta^* \subset \underline \Theta_\beta$.
It is open whether or not
\[              \hdim[\gamma(\underline \Theta_\beta)] \leq d_\beta. \]

\item   
Note that  $(0,2]  = \overline \Theta_{\beta_*}
  \cup \underline \Theta_{\beta_*}^*.$  It follows that
\[  \hdim(\gamma(0,2]) = d_{\beta_*} = d. \]
Hence, Beffara's theorem on the dimension of the path \cite{Beffara,Law1}
is a particular case of the theorem.

\item The statements about the dimension of
$\gamma(\Theta_\beta),\gamma(\overline \Theta_\beta),
\gamma(\underline \Theta_\beta^*)$  
are independent of the parametrization of the curve.

\item  Using the Markov property for $\SLE$ it is not hard to show
that with probability one, either $\Theta_\beta$ is dense in
$(0,\infty)$ or it is empty.  Also, $\hdim[\gamma(\Theta_\beta
\cap [t_1,t_2])]$ is the same for all $0 < t_1 < t_2\leq 2$.   
 In particular, in order to prove the lower bound  on
dimension, it
suffices to prove that for all $\alpha < d_\beta$, 
\[  \Prob\left\{\hdim[\gamma(\Theta_\beta \cap [1,2])] 
   \geq \alpha \right\} > 0 . \]
This is what we will do in Section \ref{lowersec}.
The proof proves the
slightly stronger (for $\kappa > 4$)  result
\[  \Prob\left\{\hdim[\Half \cap \gamma(\Theta_\beta \cap [1,2])] 
   \geq \alpha \right\} > 0 . \]

\item  If $\kappa =8$, $\beta_* = \beta_+ =1$ and $\hdim[
\gamma(\Theta_1)] = 2$.  This is related to the fact that
this is the hardest case to establish the existence of
the curve; the curve is almost surely not H\"older continuous (in the capacity parameterization) when $\kappa = 8$ \cite{JL}. For other values of $\kappa$, $\beta_* < \beta_+ < 1$. 
\end{itemize}

\subsection{Multifractal spectrum of harmonic measure}  \label{harmonsec}
Let $\Theta_\alpha^\hm$ be defined as in Section~\ref{hmsec}.
Let
\[  F_{\rm{tip}}(\alpha) := d_{1- 1/\alpha}
 =\alpha\left(1-\fr{4}{\kappa}\right) +
\fr{(4+\kappa)^2}{8 \kappa}-\fr{\kappa}{8} \left(
  \frac{\alpha^2}{2\alpha -1} \right),  
\]
and let  $\alpha_-, \alpha_*,\alpha_+$ correspond to $\beta_-,
\beta_*,  \beta_+$ through the
relation 
\[
\alpha=\frac{1}{1-\beta}.
\] 
\begin{remark}
We can compare the function $F_{\rm{tip}}$ with the conjectured almost sure bulk spectrum for $\SLE_\kappa$ given by
\[
F_{{\rm bulk}}(\alpha)=\alpha+ \frac{(4+\kappa)^2}{8 \kappa}-
  \fr{(4+\kappa)^2}{8 \kappa}\left( \frac{\alpha^2}{2\alpha -1}
\right).
\] 
\end{remark}
\begin{theorem} \label{harmtheorem}
Suppose $\alpha_{-} \le \alpha \le \alpha_+$. For chordal $\SLE_\kappa$, with probability one, 
\[  \hdim[\gamma(\Theta_\alpha^\hm)]  = F_{\rm{tip}}(\alpha). \]
\end{theorem}
\begin{proof}
This is an immediate
corollary of Theorem \ref{bigtheorem} and Lemma
\ref{aug13.lemma}.
\end{proof}
Theorem~\ref{harmtheorem} combined with
\eqref{best2} gives  gives some information on $\tilde \Theta_\alpha^\hm$. In Section \ref{newsec}, we will use
the forward Loewner flow to give a proof of the following. 
\begin{theorem}  \label{newtheorem}
If $0< \kappa < 8$ and $1/2 \le \alpha \le \alpha_*$, then with probability one
there exists a set $V$ such that  $\hdim[\gamma(V)] \le F_{\rm{tip}}(\alpha)$
and for $t \not\in V$, $\gamma(t) \in \Half$, 
\begin{equation}  \label{aug13.5}
    \tilde \mu(t,2^{-n}) \lap  2^{-n\alpha}, 
  \;\;\; n \rightarrow \infty.
\end{equation}  
\end{theorem}
Let \[
\tilde{T}^{{\rm hm}}_{\alpha}=\{t \in (0,2]: \tilde{\mu}(t, 2^{-n}) \gap 2^{-\alpha n}, \, \gamma(t) \in \mathbb{H}\}\]
and note that Theorem~\ref{newtheorem} combined with \eqref{best2} and Theorem~\ref{harmtheorem} implies that for each $\alpha_- \le \alpha < \alpha_*$ with probability one  
\[\hdim[\gamma(\tilde{T}^{{\rm hm}}_{\alpha})] = F_{\rm{tip}}(\alpha).\]
Indeed, it follows directly from \eqref{best2} and Theorem~\ref{harmtheorem} that the lower bound on the dimension holds with probability one. To get the upper bound, notice that $\tilde{T}^{{\rm hm}}_{\alpha}$ is contained in $\{t \in (0,2] : \tilde{\mu}(t,2^{-n}) \gap_{{\rm i.o.}} 2^{-n\alpha}\}$, which, for those $t$ such that $\gamma(t) \in \Half$, in turn is contained in the set $V$ from Theorem~\ref{newtheorem}.




\subsection{Parameters}  \label{parasec}

In the statement of the main theorem, $\beta$ and $\rho$
were the parameters used.  However, in deriving the result
it is useful to consider a number of other parameters.
Let
\[ r_* =  \min\left\{1,\frac 8 \kappa \right\},
    \;\;\;\;
   r_c = \frac 12 + \frac 4 \k,   \]
and note that
\[  0 < r_* \leq  r_c   , \]
where the second inequality is strict unless $\kappa = 8$.   
Let $r < r_c$;  we  define
$\lambda,\zeta,\beta,\rho$ as functions of $r$.

Let
\begin{equation}  \label{sept16.1}
  \lambda = \lambda(r)  =  r \, \left(1 + \frac \k 4\right)
   - \frac{\k r^2}{8} .
\end{equation}
 We  write $\lambda_* = \lambda(r_*)$,
and similarly for other parameters.
As $r$ increases from $-\infty$ to $r_c$, $\lambda$
increases from $-\infty$ to 
\[  \lambda_c  = 1 + \frac{3 \k}{32}
   + \frac{2}{\k} . \]
Since
the relationship is injective, we can write either
$\lambda(r)$ or $r(\lambda)$.
Solving the quadratic equation gives
 \[  r(\lambda) = \frac{4 + \k - \sqrt{(4+\k)^2 - 8\lambda \k}}
     {\k}. \]
Also,
\[\lambda(0) = 0 , \;\;\;\;\lambda_*  =
d .
\]
Let
\begin{equation}  \label{sept16.2}
  \zeta = \zeta(r) =  r - \frac{\k r^2}{8}
  = \lambda(r) -\frac{\k r}{4},
\end{equation}
and note that
\[   \zeta_* = 2-d
.\]
We can write $\zeta$ as a function of $\lambda$,
\[  \zeta(\lambda) = \lambda +
\frac{ \sqrt{(4+\k)^2 - 8\lambda \k}-4-\k}
     {4}.
  \]

We now briefly discuss some results from \cite{Law1,JL}. The
reverse-time Loewner flow $h_t$ (see Section
\ref{goodsec} for definitions) has the property that
for fixed $t$, the distribution of $|h_t'(z)|$ is
the same as that of $|\hat f_t'(z)|$. For the reverse-time
flow, if $r \in \R$ and $\lambda,\zeta$ are defined as
above, then
\[     |h_t'(z)|^\lambda \, Y_t(z)^{\zeta} \, [\sin
  \arg Z_t(z)]^{-r} \]
is a martingale.  Typically one expects $Y_t(i)
 \asymp \sqrt t$ and $ \sin
  \arg Z_t(i) \asymp 1$.  If this is true, then the
martingale property would imply
\[       \E\left[|\hat f_{t^2}'(i)|^\lambda
 \right] =      \E\left[|h_{t^2}'(i)|^\lambda\right]
 \asymp    t^{-\zeta}. \]
It turns out that this argument can be carried out
if $r < r_c$, and this is the starting point for
determining the multifractal spectrum.

We define  $\beta = \beta(r)$ by 
the relation
\[    \frac{d \zeta}{d \lambda} = - \beta . \]
A straightforward calculation gives
\[          \beta(r) = -1 + \frac{\k}{4 + \k - \k r},\;\;\,\,\,\,
  \k \,r(\beta)  = 4 + \kappa - \frac \k {\beta + 1}. \] 
Note that $\beta$ increases with $r$ with
\[ \beta({-\infty}) = -1, \;\;\;\;
  \beta(0) = - \frac{4}{4 + \k} = \beta_\#, \;\;\;\;
\beta(r_*) = \beta_*, \;\;\;\;
   \beta_c
 = 1 ,  \]
where $\beta_\#,\beta_*$ are as defined in
the previous section.  Roughly speaking, 
$\E[|\hat f_{t^2}'(i)|^\lambda]$ is carried on an event
on which $|\hat f_{t^2}'(i)| \approx t^{\beta}$ and
\begin{equation}  \label{apr10.2}
   \Prob\left\{ |\hat f_{t^2}'(i)| \approx
              t^{\beta}\right\} \approx  t^{-(\zeta
  + \lambda \beta)}. 
\end{equation}

We emphasize that the relation between $r, \lambda, \beta$ for
$-\infty < r < r_{c}$ is
bijective  and in order to specify the values of the parameters
it suffices to give the value of any one of these.  For example,
we could choose $\beta$ as the independent variable and write
$r(\beta), \lambda(\beta)$.   This is the natural approach when
proving Theorem \ref{bigtheorem}, but the formulas tend to
be somewhat simpler if we choose $r$ to be the independent
variable.

From \eqref{apr10.2}, it is natural to define
\[ \rho = \rho(r) = \zeta(r) + \lambda(r) \, \beta(r)
  = \frac{\k^2 \, r^2}{8(4 + \k - \k r)}. \]
We can also write $\rho$ as a function of $\beta$
and a computation gives \eqref{nov10.1}.
Note that
\[  \frac{d\rho}{d\beta} = \frac{d \zeta}{d\lambda}
    \, \frac {d\lambda}{d\beta} + \lambda + \beta \, \frac{d\lambda}
      {d\beta} = \lambda. \]

Let $r_+,r_-$ denote the two values of $r$ for which
$\zeta(r) + \lambda(r) \, \beta(r) = 2$ with corresponding
values $\beta_+ = \beta(r_+),\beta_- = \beta(r_-)$.   Then
\[  r_\pm =    \frac 4 \k \, \left[
    -2 \pm \sqrt{8 + \k}\right], \]
and $\beta_+,\beta_-$ are given as in \eqref{nov10.2}--\eqref{nov10.3}.
Note that if $\kappa \neq 8$, then $r_+ < r_c$.

Define
\[  d_\beta =  \frac{2 - \rho(\beta)}{1-\beta} =
\frac{2-(\zeta + \beta \lambda)}
      {1 - \beta} . \]
Note that $d_\beta$ is maximized at $\beta = \beta_*$
(interpreted as a limit for $\kappa = 8$) with
$d_* = d$.  We can also define $d$ as a function
of $r$,
\[  d(r) = 1 + \frac{\kappa - \frac{\kappa^2r^2}{8}}
   {8 + \kappa - 2 \kappa r}. \]
Straightforward differentiation shows that $d'(r) = 0$ implies
$r=1$ or $r=8/\kappa$.  Note that $1 = 8/\kappa = r_+$ if
$\kappa = 8$ and
\[          1 < r_+ < \frac{8}{\kappa}, \;\;\;\; \kappa < 8,\]
\[    \frac{8}{\kappa} < r_+ < 1 , \;\;\;\; \kappa > 8.\]
From this we can see that $d(r)$ achieves its maximum on $(-\infty,r_+)$
at $r = r_*$; in fact, $d(\beta)$ increases for $\beta < \beta_*$
and decreases for $\beta_* < \beta < \beta_+$.

In order to match the notation of \cite{Law1},
let 
\begin{equation}  \label{sept16.3}
           q = r_c - r  = \frac 12 + \frac 4 \k
 - r .
\end{equation}
Obviously, $q > 0$ if and only if $r < r_c$.  For future
reference we note that
\begin{equation}  \label{sept16.5}
 \frac{1-2q}{1+2q} = \beta.  
\end{equation}

\section{Proof of the main theorem}  \label{proofsec}

In this section we will
present the proof of the theorem relying on estimates
about moments of derivatives of the map
$\hat f$.
The upper bounds 
are proved in Section \ref{uppersec}, and 
 the lower bound is proved in Section
\ref{lowersec}.

\subsection{Upper bounds}  \label{uppersec}  

In this subsection (and this subsection only)
we  
write
\[    \hat f_{j,n} = \hat f_{(j-1)2^{-2n}}.
  \] 
For each $t \in [0,\infty),$ we associate a dyadic time by
defining
\[  t_{n} = t_n(t) = \frac{j-1}{2^{2n}} \;\;\; \mbox{ if }
 \;\; \frac{j-1}{2^{2n}} \leq t  <\frac{j}{2^{2n}} .\]
We fix $s$ with $0 < s  \leq 2$ and allow
  constants to depend on $s$.

The next theorem states
the derivative estimates that we will use
for the upper bounds; a proof can be found in \cite{JL}.

\begin{theorem} \cite{JL} \label{nov9.prop20}
If $r < r_c$, there exists $c < \infty$ such that
 for all $t \geq 1$, 
\begin{equation}  \label{upper2}
 \E\left[|\hat f_{t^2}'(i)|^{\lambda}\right]
\leq c \, t^{-\zeta}.
\end{equation}
\end{theorem}

\begin{corollary}   If $ \beta \geq \beta_\#$,
 there is a $c  < \infty$
such that if 
\[  N_{n,\beta} =   \sum_{s2^{2n} \leq j \leq
 2^{2n+1}}
      1\left\{  |\hat f_{j,n}'(i2^{-n}) |
\geq 2^{n\beta} 
   \right\},\] then
\begin{equation}  \label{nov15.2}
       \E[N_{n,\beta}] \leq c \, 2^{n  (2 - \rho)}. 
\end{equation}

\end{corollary}

\begin{proof}  
The range $\beta \geq \beta_\#$
corresponds to 
  $\lambda \geq 0$.   Hence, by Chebyshev's inequality,
\begin{eqnarray*}
   \Prob\{|\hat f_{j,n}'(i2^{-n}) | \,
 \geq 2^{n{\beta}} \}
   &  \leq&  2^{-n\beta \lambda} 
 \E\left[|\hat f_{j,n}'(i2^{-n}) |^{\lambda}\right] \\
 & = & 2^{-n\beta\lambda} \, \E\left[|\hat f_{j}'(i)|^\lambda \right] 
 \leq c   \, 2^{-n \beta \lambda}
  \,  j^{-\zeta/2}, 
\end{eqnarray*}
and hence 
\begin{eqnarray*}
  \E\left[N_{n,\beta} \right]
 & =  &  \sum_{s2^{2n}  \leq j \leq
 2^{2n+1 }}
      \Prob\{  |\hat f_{j,n}'(2^{-n})| 
 \geq 2^{n\beta}  \}  \\
 & \leq  &  \, 2^{-n \beta \lambda} \, \sum_{s2^{2n}  \leq j \leq
 2^{2n+1 }} j^{-\zeta/2}\\
 & \leq & c  \, 2^{n(2-\beta \lambda - \zeta)}
  = c  \, 2^{n(2-\rho)}
.
\end{eqnarray*}
 
\end{proof}

\begin{corollary}  If $\beta < \beta_\#$,
 there is a $c  < \infty$
such that if 
\[   N_{n,\beta}^* =  \sum_{s2^{2n } \leq j \leq
 2^{2n+1}}
      1\left\{  |\hat f_{j,n}'(i2^{-n}) |
\leq 2^{n\beta} 
   \right\},\] then
\begin{equation}  \label{nov15.2.alt}
       \E[N_{n,\beta}^*] \leq c \, 2^{n  (2 - \rho)}. 
\end{equation}

\end{corollary}

\begin{proof}  This is proved in the same way using
$\lambda < 0$. 
\end{proof}

The standard technique to find
upper bounds for Hausdorff dimension uses
 an appropriate sequences of covers for a set.
We will now describe the covers that we will use.   
 Let
\[  I(j,n) = \left[\frac{j-1}{2^{2n}}, \frac j{2^{2n}}\right].\]
If $b,\overline b \in \R$ with
$-1 < b < \overline b < 1$, 
let $\overline
\ball (j,n,\bar b)$ be the closed disk in $\C$
 of radius 
$2^{n(\bar b-  1)}$ centered at  $\hat f_{j,n}(i2^{-n})$, and
let
\[     I_n(s,b) = \bigcup I(j,n) \;\;\;\;
   \ball_n(s,b,\overline b) = \bigcup \overline \ball(j,n,
\overline b) , \] 
where in each case the union is over   $s2^{2n} \leq
j \leq 2^{2n+1}$ with
$
   |\hat f_{j,n}'(i2^{-n})| \geq 2^{nb} $.  Let
\[    I^m(s,b)  = \bigcup_{n=m}^\infty I_n(s,b) ,\;\;\;\;
     \ball^m(s,b,\overline b) 
 = \bigcup_{n=m}^\infty \ball_n(s,b,\overline b)  . \]

\begin{lemma}  \label{apr8.lem1}
If $-1< b < \beta < b_1 < \overline b < 1$, then
for each $m$, 
\[    \overline \Theta_\beta \cap (s,2]
 \subset  
    I^m(s,b), \]
\[    \gamma\left(\overline \Theta_\beta
\cap 
  \tilde  \Theta_{b_1} \cap (s,2]\right)
   \subset \ball^m(s,b,\overline b). \]
\end{lemma}

\begin{proof}
Suppose $t \in \overline 
\Theta_\beta \cap (s,2]$. By Proposition \ref{distortprop2},
 there exists a subsequence
$n_j \rightarrow \infty$   such that
\[ |f_{t_{n_{j}}}'(i2^{-n_j})| \geq  2^{n_jb}  . \]
In other words, there is a sequence $n_j$ such that
$I(t_{n_j} 2^{2n_j},n_j) \in I_{n_j}(s,b)$. 
This proves the first assertion. 

If $t \in \overline \Theta_\beta \cap \tilde
  \Theta_{b_1}   \cap (s,2]$ and
 $b_1< u < \overline b$, then \eqref{110} shows
that for all $n$ sufficiently
large,
\[   |\gamma(t) - \hat f_{t}(i2^{-n})|
  \leq  2^{(u-1) n}.\]
The triangle inequality gives
\[  |\gamma(t) - \hat f_{t_n}(i2^{-n})| \leq 
|\gamma(t) - \hat f_{t}(i2^{-n})|
 + |\hat f_t(i2^{-n}) -
  \hat f_{t_n}(i2^{-n})|, \]
and estimating as in \eqref{nov20.20}, we have
for $n$ sufficiently large
\[ |\hat f_t(i2^{-n}) -
  \hat f_{t_n}(i2^{-n})| \leq  2^{(u-1) n}.\]
Hence, for $n$ sufficiently large,
\[    |\gamma(t) - \hat f_{t_n}(i2^{-n})|  
\leq  2^{n_j (\bar b-1)} .\]
This implies for all $j$ sufficiently large,
\[    \gamma(t) \in \ball(t_{n_j}2^{2n_j},b,\overline b).\]
\end{proof}

\begin{proposition}   \label{prop70} If $
\beta \geq \beta_\#$, then  with
probability one, 
\begin{equation} \label{101}
   \hdim(\overline \Theta_\beta \cap
(s,2]) \leq \hat d_\beta, 
\end{equation}
Moreover, if $\beta > \beta_+$,
then with probability one
\[   \overline \Theta_\beta \cap
(s,2] = \emptyset.\]
\end{proposition}

\begin{proof}
It suffices
to consider $  \beta_\#< \beta < 1$. Suppose $
\beta_\# < b <\beta < 1$. 
Using the cover from Lemma \ref{apr8.lem1}, 
we  get
\[  \haus^\alpha \left[\overline \Theta_\beta \cap (s,2]\right] 
\leq \sum_{n=m}^\infty
              N_{n,b} \, 2^{-2\alpha n}, \]
and hence \eqref{nov15.2} implies 
\[  \E \left(\haus^\alpha \left[\overline \Theta_\beta \cap (s,2]\right] 
  \right) 
\leq c \sum_{n=m}^\infty
               2^{n(2-\rho(b))}
   \, 2^{-2\alpha n}. \]
The sum goes to zero, provided that  $2 \alpha > 2 -\rho(b)$,
and hence with probability one
\[           \haus^\alpha(\overline \Theta_\beta \cap
(s,2]) = 0 , \;\;\;\;\;  \alpha > \frac{2 -\rho(b)}2.\]
Letting $b \rightarrow \beta$ gives \eqref{101}.

For the second assertion, note that
\[ \Prob\left\{\overline \Theta_\beta \cap (s,2] \neq
\emptyset\right\} \leq \sum_{n=m}^\infty  \E[N_{n,b}]
         \leq c\sum_{n=m}^\infty 2^{n(2-\rho(b))}. \]
If $\beta > \beta_+$, then $\rho(\beta) > 2$
and we can find $b < \beta$ with
$\rho(b) > 2$.
\end{proof}

\begin{lemma} \label{prop71} If $\beta_\# \leq   \beta < b_1 < 1$,
then with probability one, 
\begin{equation}  \label{103}
    \hdim\left[\gamma(\overline \Theta_\beta
\cap \tilde \Theta_{b_1}  \cap
(s,2])\right] \leq  \frac{2-\rho(\beta)}{1-b_1}. 
\end{equation}
\end{lemma}

\begin{proof}
Choose $b,\overline b$ with $\beta_\# < b
 < \beta < b_1 < \overline b < 1.$
Using the cover from Lemma \ref{apr8.lem1}, 
we  get
\[ 
 \haus^\alpha\left[\gamma(\overline \Theta_\beta
\cap \tilde \Theta_{b_1}  \cap
(s,2])\right] \leq
            \sum_{n=m}^\infty  N_{n,b} \, 2^{\alpha (\bar b-1)n}, \]
and hence \eqref{nov15.2} implies 
\[ 
 \E\left(\haus^\alpha\left[\gamma(\overline \Theta_\beta
\cap \tilde \Theta_{b_1}  \cap
(s,2])\right]\right) \leq
            c\sum_{n=m}^\infty 2^{n(2-\rho(b))}
 \, 2^{\alpha (\overline b-1)n}. \]
The sums on the right go to zero provided that
\[   
    \alpha > \frac{2 - \rho(b)}{ 1-\overline b} , \]
respectively.  
We now choose a sequence of values for $b,\overline b$
that converge to $\beta,b_1$ to conclude \eqref{103}.
\end{proof}

\begin{proposition}  If $\beta_* \leq \beta \leq 1$, then
with probability one, 
\begin{equation}  \label{102}
       \hdim\left[\gamma(\overline \Theta_\beta
   \cap (s,2])\right] \leq d_\beta . 
\end{equation}
\end{proposition}

\begin{proof}
 
If $\beta =b_0
 < b_1 < b_2 < \ldots < b_k < 1$ with
$b_k > b_+$, then 
\[ \overline \Theta_\beta    = \bigcup_{j=1}^k
            (\overline \Theta_{b_{j-1}} \cap \tilde
   \Theta_{b_j}) .\] 
Therefore, \eqref{103} implies
\[  \hdim\left[\gamma(\overline \Theta_\beta
   \cap (s,2])\right]  \leq \max\left\{\frac{2-\rho(b_{j-1})}
    {1 - b_j}:j=1,\ldots,k\right\}. \]
By taking finer and finer partitions and using
the continuity of $\rho$, we see that
\[  \hdim\left[\gamma(\overline \Theta_\beta
   \cap (s,2])\right]  \leq \sup_{b \geq \beta}
             \frac{2-\rho(b)}{1-b} = d_\beta. \]
The last equality uses $\beta \geq \beta_*$ and the
fact which can easily be verified (see Section
\ref{parasec}) that the function
\[            b \mapsto    \frac{2-\rho(b)}{1-b} , \]
is decreasing on the interval $[\beta_*,\beta_+]$.
\end{proof}

For $\beta < \beta_\#$ we use a slightly different
cover.  Let $I(j,n)$ be as above and 
\[     I_n^*(s,b) = \bigcup I(j,n) \;\;\;\;
   \ball_n^*(s,b,\overline b) = \bigcup \ball(j,n,
\overline b) , \] 
where in each case the union is over   $s2^{2n} \leq
j \leq 2^{2n+1}$ with
$
   |\hat f_{j,n}'(i2^{-n})| \leq 2^{nb} $.  Let
\[    I^m_*(s,b)  = \bigcup_{n=m}^\infty I_n^*(s,b) ,\;\;\;\;
     \ball^m_*(s,b,\overline b) 
 = \bigcup_{n=m}^\infty \ball_n^*(s,b,\overline b)  . \]

\begin{lemma}
If $-1 < \beta < b < \overline b,$
\[    \underline \Theta_\beta \cap (s,2]
  \subset I^m_*(s,b) , \]
\[   \gamma\left(\underline \Theta_\beta^*  \cap (s,2]\right)
   \subset \ball_*^m(s,b,\overline b). \]
\end{lemma}

\begin{proof}
Suppose $t \in \underline 
\Theta_\beta \cap (s,2]$. By Proposition \ref{distortprop2},
 there exists a subsequence
$n_j \rightarrow \infty$   such that
\begin{equation}  \label{apr9.1}
 |\hat{f}_{t_{n_{j}}}'(i2^{-n_j})| \leq  2^{n_jb}  . 
\end{equation}
In other words, there is a sequence $n_j$ such that
$I(t_{n_j} 2^{2n_j},n_j) \in I_{n_j}^*(s,b)$. 
This proves the first assertion. 

Now suppose $t \in \underline \Theta_\beta^* \cap (s,2]
 \subset \underline \Theta_\beta \cap (s,2]$.  Then
there exists a sequence $n_j$ such that both \eqref{apr9.1}
holds and
\[      v_{t}(2^{-n_j}) \leq 2^{n_j(b-1)}. \]
Using the triangle inequality as in Lemma \ref{apr8.lem1},
we see that
\[|\gamma(t) - \hat f_{t_{n_j}}(i2^{-n_j})|
    \leq  v_{t}(2^{-n_j})  +
      \left|\hat f_t(i2^{-n_j}) -
\hat f_{t_{n_j}}(i2^{-n_j})\right|, \] 
and arguing as before we see that for $j$ sufficiently
large
\[ |\gamma(t) - \hat f_{t_{n_j}}(i2^{-n_j})|
    \leq  2^{n_j(\overline b -1)}, \]
and
\[ \gamma(t) \in \ball^*(t_{n_j}2^{2n},b,\overline b). \]
\end{proof}

\begin{proposition}  If $\beta < \beta_\#$, then with
probability one,
\[         \hdim(\underline \Theta_\beta \cap (s,2])
               \leq \hat d_\beta,\]
\[   \hdim(\gamma\left(\underline \Theta_\beta^* \cap (s,2]\right))
               \leq d_\beta.\]
Moreover, if $\beta < \beta_-$, with probability one,
\[           \underline \Theta_\beta \cap (s,2]
  = \emptyset.\]
\end{proposition}

\begin{proof} This is proved in the same way as
Proposition \ref{prop70} and Lemma \ref{prop71}
using \eqref{nov15.2.alt}.
\end{proof}

\subsection{Lower bound}  \label{lowersec}

In this subsection we prove the lower bound for 
the dimension in \eqref{sept18.0}.
We fix $r$ such that $\rho = \lambda \beta + \zeta < 2$
and recall that $r < r_c$.
As has been pointed out, it suffices to show that with
positive probability
\begin{equation}  \label{nov15.3}
  \hdim(\Theta_\beta \cap [1,2]) \geq
      \frac{2-\rho}{2} , \;\;\;\;\;
   \hdim[\gamma(\Theta_\beta \cap [1,2])]
  \geq \frac{2-\rho}{1-\beta}. 
\end{equation}
We will use a standard technique of Frostman to show
that with positive probability there exist
nontrivial positive measures $\mu,\nu$ whose supports
are contained in  $\Theta_\beta \cap [1,2]$
and $\gamma(\Theta_\beta \cap [1,2])$, respectively, such
that
\[  \energy_\alpha(\mu) < \infty , \;\;  \alpha < 
   \frac{2-\rho}{2};\;\;\;\;\; 
  \energy_\alpha(\nu) < \infty, \;\;  \alpha
   < \frac{2-\rho}{1- \beta}, \]
where
\[ \energy_\alpha(\mu) = \int \int \frac{\mu(dx) \, \mu(dy)}
   {|x-y|^\alpha} \]
is the energy integral.
It is well known that this implies \eqref{nov15.3}.
For this subsection, we will adopt a different
notation than in the previous subsection.
  We let
\[    \hat f_{j,n} = 1 + \frac{(j-1)}{n^2}, \;\;\;\;
  j=1,2,\ldots,n^2 .  \]
We will be studying $|\hat f_{j,n}'(i/n)|$.
The proof considers a subset of
times in $\Theta_\beta \cap [1,2]$ that behave in
some sense nicely.
The hard work is Theorem \ref{theorem.nov7} which will
be proved in Section \ref{estsec}.  
This  theorem discusses the existence of some events
 $E_{j,n}$ on which 
 \[
 |\hat f_{j,n}'(iy)| \approx y^{-\beta}, \;\;\;
  n^{-1} \leq y \leq 1.  \] 
The definition of the events (``good times'') will be left
for Section \ref{estsec}.

\begin{theorem}  \label{theorem.nov7}
Suppose  $
\rho = \lambda \beta + \zeta < 2$.
There exist $c < \infty$, a subpower function $\psi$, and events
\[   E_{j,n}, \;\;\;n=1,2,\ldots, \;\; j=1,\ldots,n^2, \]
such that the following hold.  Let $E(j,n) = 1_{E_{j,n}}$ and 
\[     F(j,n) = n^{\zeta-2} \, |\hat f_{j,n}'(i/n)|^\lambda \, 
  E(j,n).\]
\begin{itemize}
\item If $1 \leq j \leq n^2$, then on the event $E_{j,n}$,
\begin{equation}  \label{sept17.1}
  \psi(1/y)^{-1} \, y^{-\beta} \leq
   |\hat f_{j,n}'(iy)| \leq  \psi(1/y)\,
 y^{-\beta} , \;\;\;\;
          n^{-1} \leq y \leq 1. 
\end{equation}
 \item If $1 \leq j \leq k \leq n^2$, 
\begin{equation}  \label{sept17.1.1}
    c_1 \, n^{-2} \leq  \E\left[F(j,n) \right]
     \leq n^{\zeta - 2} \, \E\left[ |\hat f_{j,n}'(i/n)|^\lambda 
 \right] \leq c_2 \, n^{-2}, 
\end{equation}
\item If $1 \leq j \leq k \leq n^2$,
\begin{equation}  \label{sept17.1.2}
  \E[F(j,n) \, F(k,n)] \leq  n^{-4} \, \left(\frac{n^2}{k-j+1}\right)
     ^{\frac{\rho}{2}} \, \psi\left(\frac{n^2}{k-j+1}\right).
\end{equation}
\item If $1 \leq j < k \leq n^2$ and $E(j,n) \, E(k,n) = 1$,
\begin{equation}  \label{sept17.1.3}
  \left|\hat f_{j,n}(i/n) - \hat f_{k,n}(i/n)\right|^{-\frac{2-\rho}{1-\beta}}
   \leq 
  \left(\frac{n^2}{k-j+1}\right)
     ^{\frac{2-\rho}{2}} \, \psi\left(\frac{n^2}{k-j+1}\right)
^{d_\beta}.
\end{equation}
\end{itemize}
\end{theorem}

\begin{proof}  This theorem combines Propositions \ref{nov9.10}
and \ref{nov9.10.1} proved in Section \ref{estsec}. 
\end{proof}

\begin{proposition} Under the assumptions of Theorem \ref{theorem.nov7},
with positive probability there exists $A \subset [1,2]$ such that
for $t \in A$,
  \begin{equation}  \label{nov7.2}
 \frac{1}4 \, y^{-\beta} \, \psi(1/y)^{-1} \leq
   |\hat f_t'(iy)| \leq 4 \, y^{-\beta} \, \psi(1/y) , \;\;\;\;
   0 < y \leq 1.
\end{equation}
and such that
\[  \hdim(A) \geq \frac{2-\rho}{2},\;\;\;\;\;
  \hdim[\gamma(A)] \geq \frac{2-\rho}{1-\beta}. \]
\end{proposition}

\begin{proof}[Proof assuming Theorem \ref{theorem.nov7}.]
We use a now standard argument to show that with positive
probability a ``Frostman measure'' of appropriate dimension
can be put on the set of $t$ satisfying \eqref{nov7.2}.
The proof is very similar to that of \cite[Lemma 10.3]{Law1}
so we omit some of the details.

Let $\mu_{j,n}$ denote the random measure on $\R$ that is a multiple
of Lebesgue measure on $I(j,n) := [1+(j-1)n^{-2},1+jn^{-2}]$
where the multiple is chosen to that  $\|\mu_{j,n}\| = F(j,n)$.
Let $\nu_{j,n}$ denote the random measure on $\C$ that is a multiple of
Lebesgue measure on the disk of radius $n^{\beta - 1} \,\psi(n^2)^{-1}/
4$ centered at
 $\hat f_{j,n}(i/n)$ where the constant is chosen so
that $\|\nu_{j,n}\| = F(j,n)$.  Let $\mu_{n} = \sum_{j=1}^{n^2}
\mu_{j,n}, \nu_n = \sum_{j=1}^{n^2}
  \nu_{j,n}$.   Note that
\[  \|\mu_n\| = \|\nu_n\| = \sum_{j=1}^{n^2} F(j,n). \]

From \eqref{sept17.1.1} and \eqref{sept17.1.2}, we see that 
\[  \E\left[\|\mu_n\|\right] \geq c_1, \;\;\;\;
    \E\left[\|\mu_n\|^2\right] \leq c_2. \]
Hence 
\[
\Prob \left\{ \| \mu_n\| > 0 \right\} \ge c >0,
\]
uniformly in $n$.
From \eqref{sept17.1.2} and \eqref{sept17.1.3}, we can show
that there is a $C_\alpha$ such that

\[          \E[\energy_\alpha(\mu_n)] \leq C_\alpha ,\;\;\;\;\;
   \alpha < \frac{2-\rho}{2}, \]
\[         \E[\energy_\alpha(\nu_n)] \leq C_\alpha, \;\;\;\;\;
     \alpha < \frac{2-\rho}{1-\beta}. \]

We let $\mu$ denote a subsequential limit of the $\mu_n$ which with
positive probability we know is nontrivial and satisfies
$\energy_\alpha(\mu) < \infty$ for all $\alpha < (2-\rho)/2$.
Hence,
\[   \hdim[\supp \mu] \geq \frac{2-\rho}{2}.\]
Similarly, let $\nu$ denote a subsequential limit of the $\nu_n$ 
which is nonzero with positive probability and satisfies
\[   \hdim[\supp \nu] \geq \frac{2-\rho}{1-\beta} . \]

We claim that every $t \in \supp \mu$ satisfies \eqref{nov7.2}. 
 Indeed, the construction shows that
if $t \in \supp \mu$, then 
 there is a subsequence $n_k \rightarrow \infty$
and $j_k \in \{1,\ldots,n_k^2\}$ such that
$E(j_k,n_k) = 1$ and 
\begin{equation}  \label{nov7.1}
           \lim_{k \rightarrow \infty}
        \frac{j_k-1}{n_k^2} = t . 
\end{equation}
Suppose for some $t\in [1,2]$ and $0 < y \leq 1$, we had
\[             |\hat f_t'(iy)| \geq  4 \, y^{-\beta} \, \psi(1/y). \]
Continuity would imply that for all $s$ in a neighborhood of
$t$,
\[        |\hat f_s'(iy)| \geq 2 \, y^{-\beta} \, \psi(1/y). \]
This implies that there is no sequence $(j_k,n_k)$ as above with
$E(j_k,n_k)$ satisfying \eqref{nov7.1}.  A similar argument
shows that there cannot exist $t \in \supp \mu$ and
$y$ with 
  $ |\hat f_t'(iy)| \leq  (1/4) \, y^{-\beta} \, \psi(1/y)^{-1}$, 
 and this gives \eqref{nov7.2}.
Similarly, $\supp \nu$ is contained in $\gamma(A')$, where $A'$
denotes the set of $t\in [1,2]$ satisfying \eqref{nov7.2}. 
\end{proof}

\section{Estimating the moments}  \label{estsec}

In this section  $r < r_c$  with corresponding values
of $\zeta,\beta,\rho,\lambda$.  All constants may depend
on $r$. Let us give an overview of the section. We begin by discussing the reverse-time Loewner flow and how it relates to $f_t$. We then go on to define the ``good times'' which, roughly speaking, are $T$ for which the reverse flow driven by $t \mapsto V_{T-t}-V_t$ behaves in some sense nicely. (Here, $V_t$ is a two-sided Brownian motion.) We make this precise in a number of lemmas that show how $|h'_t|=|h'_{T,t}|$ can be controlled on the event that $T$ is ``good''. The needed correlation estimates can then be derived using moment bounds from \cite{Law1, JL}.

\subsection{Reverse Loewner flow}

Here we state the basic lemma that relates the reverse 
Loewner flow to the forward flow for $\SLE$.    We will
estimate the moments for $h,\tilde h$ rather than for
$\hat f$.

If $V_t$ is a continuous function,
 define $g_t$ to be the solution to the
forward-time (chordal) Loewner equation
\begin{equation}  \label{forward.1}
        \p_t g_t(z) = \frac{a}{g_t(z) - V_t} ,\;\;\;\;
     g_0(z) = z.
\end{equation}
Let $f_t(z) = g_t^{-1}(z), \hat f_t(z) =
f_t(z + V_t)$.  If $U_t$ is another continuous
function, let $h_t$ be the solution to
the reverse-time (chordal) Loewner equation
\begin{equation}  \label{reverse.1}
     \p_t h_t(z) = \frac{a}{U_t - h_t(z)} , \;\;\; h_0
 (z) = z.
\end{equation}  The next 
lemma  relates the forward-time and reverse-time
equations; although versions of this have appeared before
we give a short proof.  We point out that this is a 
fact about the Loewner equation itself; no assumptions are
made about the function $V$ other than continuity.

\begin{lemma}  \label{loewnerlemma}
  Assume $V_t, -\infty <  t < \infty$, is a
continuous function with $V_0 = 0$.  For each
$T\geq 0 $,  let
\[        U_{t,T} =V_{T-t} - V_T.\]
Let $g_t, 0 \leq t < \infty$, be the solution to the forward-time Loewner equation \eqref{forward.1},
and let $f_t,\hat f_t$ be as above.
Let $h_{t,T}, 0 \le t < \infty$, be the solution to the reverse-time Loewner equation \eqref{reverse.1}
with $U_t = U_{t,T}$.  Let
\[   Z_{t,T}(z) = h_{t,T}(z) - U_{t,T}.\]
Then
\[   h_{T,T}(z) = \hat f_T(z) - V_T
  = \hat f_T(z) + U_T, \]
and if $0 \leq S \leq T$, $z,w \in \Half$,
\begin{equation}  \label{nov8.1}
    h_{T,T}(z) =  h_{S,S} \left(Z_{T-S,T}(z) 
\right) +U_{T-S,T} ,
\end{equation}
\begin{equation}  \label{nov8.3}
   \hat f_T(z) - \hat f_S(w) =
      h_{S,S} \left(Z_{T-S,T}(z) 
 \right) 
   - h_{S,S}(w).
\end{equation}
In particular,
\[     h_{T,T}'(z) = h_{S,S}'\left(Z_{T-S,T}(z) 
 \right) \; h_{T-S,T}'(z).\]
\end{lemma}

\begin{proof}  Fix $T$ and let $U_t = U_{t,T},
h_t = h_{t,T}$.  For $0 \leq s \leq T$, we have
\[   U_{t,s} = V_{s-t} - V_s = U_{T-s+t} - U_{T-s}.\]
 Let
$u_t(z)  = h_{T-t}(z) + V_T$.  Then $u_t$ satisfies
\[  \p_t u_t(z) = - \p_t h_{T-t}(z) = 
\frac{a}{h_{T-t}(z)- [V_{t} - V_T]}
  = \frac{a}{u_t(z) - V_t}. \]
Also, $u_T(z) = h_0(z) + V_T = z + V_T$.  By comparison,
with \eqref{forward.1}, 
\[    u_0(z) =  g_t^{-1} (z + V_T) = \hat f_t(z) , \]
and hence 
\[ h_T(z) = u_0(z) - V_T = \hat f_T(z) - V_T
  = \hat f_T(z) + U_T. \]
 
For $0 \leq s \leq T$, define $h^{(s)}_t$ by
\[                h_{s+t} = h^{(s)}_t \circ h_s. \]
By \eqref{reverse.1} we see that 
\[    \p_t [h^{(s)}_t(z)-U_s]
 = \frac{a}{U_{t+s} - h^{(s)}_t(z)}
  = \frac{a}{U_{t+s} - U_s- [h^{(s)}_t(z)-
U_s]} , \]\[
     h^{(s)}_0(z) -U_s = z - U_s. \]
Therefore,
\begin{equation}  \label{nov8.40}
     h_t^{(s)}(z) - U_s = h_{t,T-s}(z-U_s), 
\end{equation}
which implies
\begin{equation}  \label{nov8.4}
   h_{s+t}(z) = h_{t,T-s}(Z_{s,T}(z)) + U_s. 
\end{equation}
Setting $s=T-S$ and $t=S$ gives \eqref{nov8.1} and
setting $s=S,t=T-S$ gives
\[   \hat f_{T-S}(w) = 
 h_{T-S}(w) - U_{T-S,T-S} =
      h_{T-S,T-S}(w) - [U_T - U_S], \]
\[  \hat f_{T}(w) =
h_{T}(z) - U_T =  h_{T-S,T-S}(Z_{S,T}(z))
  + U_S - U_T . \]
Subtracting these equations gives \eqref{nov8.3}.
The final assertion follows from \eqref{nov8.1} and the
chain rule.
\end{proof}

The preceding
 lemma holds for all continuous $V_t$.  If $V_t$
is a standard Brownian motion, then so is $U_{t,T}$ for 
each $T$.  We get the following corollary.

\begin{lemma}
Suppose $0 < S < T$ and $g_t, \, 0 \leq t \leq T,$ is the solution to
\eqref{forward.1} where $V_t$ is a standard Brownian motion.
Suppose $U_t$ is a standard Brownian motion and $h,\tilde h$
are the solutions to
\begin{equation}  \label{reverse}
     \p_t h_t(z) = \frac{a}{U_t - h_t(z)} , \;\;\; h_0
 (z) = z, 
\end{equation}
\[     \p_t \tilde h_t(z) = \frac{a}{\tilde U_t - \tilde
 h_t(z)} , \;\;\; \tilde  h_0
 (z) = z, \]
where $\tilde U_t = U_{T-S+t} - U_{T-S}. $
Then
\[   h_{T-S + t}(z) = \tilde h_t(
    h_{T-S}(z) - U_{T-S}) + U_{T-S}. \]
Moreover, the joint distribution of the functions
\[            \left(\;\; \hat f_S'(w)\; ,\; \Im \hat f_S(w),  \;
 \hat f_T'(z)\;,\; \Im \hat f_T(w) \;, \;
 \hat f_T(z) - \hat f_S(w) \;\; \right) \]
is the same as the joint distribution of 
\[          \left(\;\;
\tilde h_S'(w)\; , \; \Im \, \tilde  h_S(w) \;,
\;  \tilde h_S'(Z)
 \; h_{T-S}'(z)\; ,  
\;\Im  \, \tilde
h_S(Z)
\;,\;  \tilde h_S(Z)
   - \tilde h_S(w)  \;\; \right), \]
where $Z = Z_{T-S}(z)  = h_{T-S}(z)- U_{T-S}$.
\end{lemma}

\subsection{Good times}  \label{goodsec}
Suppose that $T>0$   and $h_t = h_{t,T}$ 
is defined as in the proof of Lemma \ref{loewnerlemma}.
 More specifically,
$g_t$ is the solution to the forward-time Loewner equation \eqref{forward.1} with a (two-sided) Brownian
motion $V_t$ as driving function, and 
$h_t$ is the solution to the reverse-time Loewner equation \eqref{reverse.1} with
$U_t = U_{t,T}= V_{T-t} - V_T$ as driving function.  Let
\[   Z_t(z) = X_t(z) + Y_t(z) = h_t(z) - U_t. \]
 Recall from Lemma~\ref{loewnerlemma} that we have
\[h_{t+s}(z) = h_{t+s,T}(z)=h_{t,T-s}(Z_s(z)) +U_s;\]
Note that
\begin{equation}  \label{apr7.1}
   h_{t+s}'(z) =  h_{t,T-s}'(Z_s(z))
    \, h_s'(z) . 
\end{equation}
If $\psi$ is a subpower function and
$0 < \delta \leq 1$, we let
\[  \hat \psi_\delta(t) =
  \min \{\psi(t/\delta), \psi(1/t)\}
  = \left\{\begin{array}{ll} \psi(t/\delta) &
 \mbox{ if } t \leq \sqrt \delta  \\
     \psi(1/t) & \mbox{ if } t > \sqrt \delta
  \end{array} \right. \; . \]
Note that  for every subpower function $\psi$ and
every $c < \infty$, there is an $M < \infty$ such that for
all $\delta > 0$,
\begin{equation}  \label{nov8.5}
      \hat \psi_\delta(t) \leq M  \;\;\;\;
  \mbox{ if } \delta \leq t \leq c \delta \mbox{ or }
  t \geq 1/c . 
\end{equation}
Roughly speaking, $\hat \psi_\delta(t)$ is $O(1)$ for
$t$ comparable to $\delta$ or comparable to $1$ but can
be larger for other $\delta < t < 1$.

\begin{definition}
 We call 
  a time $T$ $\psi$-{\em good at} $\delta$
 if the
following five conditions hold for $h_t = h_{t,T}$
 with $\hat \psi =
\hat \psi_\delta$ and 
$ Z_t = X_t + i Y_t = Z_t(\delta i).$
\begin{equation}  \label{good.1}
  Y_{t^2} \geq t \, \hat \psi(t)^{-1},\;\;\;\;
         \delta \leq t \leq 2, 
\end{equation}
\begin{equation}  \label{good.2}
  |X_{t^2}|  \leq 
   [t+\delta]\, \hat \psi(t), \;\;\;\; 0 \leq t
  \leq 2 , \end{equation}
\begin{equation}  \label{good.3}
    (t/\delta)
 ^{\beta} \, \hat \psi(t)^{-1} \leq
    |h_{t^2}'(i\delta)| \leq (t/\delta)^\beta \,
  \hat \psi(t), \;\;\;\; \delta \leq t \leq 2 ,
\end{equation}
\begin{equation}  \label{good.4}
   t^{-\beta}  \, \hat \psi(t)^{-1}
   \leq \frac{ |h_{4}'(i\delta)|}
 {|h_{t^2}'(i\delta)| } = |h_{4-t^2,T-
t^2}'(Z_{t^2})|  \leq t^{-\beta}
   \hat \psi(t) , \;\;\;\; \delta \leq t \leq 2, 
\end{equation}

\end{definition}

This definition depends on $\psi$ and $\delta$. Note
that if $T$ is $\psi$-good at $\delta$ and $\phi$
is a subpower function with $\psi \leq \phi$, then
$\psi_\delta \leq \phi_\delta$ and 
$T$ is $\phi$-good at $\delta$.  In the
remainder of this subsection, we derive some properties
of $\psi$-good times. These will be used in the later
subsections to estimate first and second moments
for $|h'_{t^2}(\delta i)|^\lambda$ on the event that $T$
is $\psi$-good at $\delta$.

\begin{proposition} For every subpower function
$\psi$ there is a subpower function $\phi$ such that
for all $\delta > 0$, if $T$ is $\psi$-good at
$\delta$, then
\begin{equation}  \label{good.5}
|U_{t^2}| \leq  t \,  \phi(1/t), \;\;\;\;
  \delta \leq t \leq 2.
\end{equation}
\end{proposition}

\begin{proof} Let $X_t = X_t(i\delta), Y_t = Y_t(
i\delta)$.  We let $\phi$ denote a subpower function
whose value may change from line to line.
 From the
Loewner equation, we know that
\[    dX_t =  -\frac{a X_t}{X_t^2 + Y_t^2} \, dt 
                - dU_t . \]
Hence,
\[
|U_{t^2}| \leq |X_{t^2}| + a\int_0^{t^2}
  \frac{|X_s|\, ds}{X_s^2 + Y_s^2} . \]
By \eqref{good.2}, it suffices to show that
\[   \int_0^{t^2}
  \frac{|X_s|\, ds}{X_s^2 + Y_s^2} \leq  
  \phi(1/t) . \]
Using \eqref{good.1} and \eqref{good.2}, we have
\[    \frac{|X_s|}{X_s^2 + Y_s^2}  \leq \frac{\phi(1/s) }
    {\sqrt s}, \]
and hence
\[ \int_0^{t^2}
  \frac{|X_s|\, ds}{X_s^2 + Y_s^2}  \leq \int_0^{t^2}
    \frac{\phi(1/s) \, ds}{\sqrt s} =
   \int_{t^{-2}}^\infty \phi(x) \, x^{-3/2} \,dx
  = \tilde \phi(1/t)  \, t,\]
where
\[ \tilde \phi(1/t) =    t^{-1} \, \int_0^{t^2}
    \frac{\phi(1/s) \, ds}{\sqrt s}.\]
It is easy to check that $\tilde \phi$ is continuous
and decays faster than
$x^\epsilon$ for each $\epsilon$.  
\end{proof}

\begin{lemma}  If $\psi$ is a subpower function, there is
a $c > 0$ such that for every $0 < \delta \leq 1$, if
$T$ is $\psi$-good at $\delta$ and $Y_t = Y_t(\delta i)$, then 
  \begin{equation}  \label{nov7.11}
    Y_{\delta^2}  
 \geq (1+c) \, \delta  .
\end{equation}
\end{lemma}

\begin{proof}
Using \eqref{good.1}, \eqref{good.2} and the
fact that $Y_{t^2}$ increases with $t$,
 we see there
is a $c_1 < \infty$ such that
\[      |X_{t^2}
  | \leq c_1 \, \delta \leq
  c_1 \, Y_{t^2}, \;\;\;\; 0 \leq t \leq \delta.\]
The Loewner equation \eqref{reverse.1} implies
that
\[      \p_s Y_s = \frac{a Y_s}{X_s^2 + Y_s^2}
  \geq \frac{a}{c_1^2 + 1}\, \frac{1}{
  Y_s} , \;\;\; s \leq \delta^2
 , \]
from which \eqref{nov7.11} follows.
\end{proof}

\begin{lemma} \label{nov9.lemma4}
 For every subpower function
$\psi$, there is a $c$ such that
if $0 < \delta \leq 1$
and   
  $T$ is $\psi$-good at $\delta$, then
\[  \left| h_{4-t^2,
T-t^2}(Z_{t^2}) -  h_{4-t^2,
T-t^2}( \delta i) 
\right|  \geq c \, \hat \psi_\delta(t)^{-2} \, t^{1-\beta},
\;\;\;\; \delta \leq t \leq 2. \]
\end{lemma}

\begin{proof} Using \eqref{good.1} and  \eqref{nov7.11},
 we see that
there 
 is a $c_1 > 0$ such that if $\ball$
denotes the 
 open disk of radius
$c_1\, t \, \hat \psi_\delta(t)^{-1}$ about
$Z_{t^2}$, then $\delta i \not \in \ball$.
Using \eqref{good.4} and   
  the Koebe $(1/4)$-theorem we see that 
$h_{4-t^2,T-
t^2}(\ball)$ contains a disk of radius   
$ (c_1/4) \, \hat \psi_\delta(t)^{-2} \, t^{1-\beta}$ about
$ h_{4-t^2, T-
t^2}(Z_{t^2}) $.  Since $h_{4-t^2,
T-t^2}(\delta i) \not\in \ball,$ the result follows.
\end{proof}


\begin{lemma}  For all subpower functions $\psi,\phi,
 $ there
is a subpower function $\overline \psi$ such that 
if $T$ is $\psi$-good at $\delta$, then
 the following
holds for all $\delta \leq t \leq  2$.
  Suppose
\begin{equation}  \label{sept28.1}
         t \,  \phi(1/t)^{-1} \leq
  y \leq t \,  \phi(1/t) , 
\end{equation}
\begin{equation}  \label{sept28.2}
       |x|  \leq [t + \delta] \, 
\phi(1/t).
\end{equation}
Then
\[     t^{-\beta} \, \overline \psi(1/t)^{-1}\leq 
  \left|h_{4-t^2,T-t^2}'(x+iy) \right|
 \leq t^{-\beta} \, \overline \psi(1/t) .\]
\end{lemma}

\begin{proof}  
By Proposition \ref{prop34} and \eqref{good.1} and
\eqref{good.2}, 
\[      \tilde \psi(1/t)^{-1} |
   h_{4-t^2,T-t^2}'(Z_{t^2})| \leq 
  \left|h_{4-t^2,T-t^2}'(x+iy) \right|
 \leq  \tilde  \psi(1/t)\, |
   h_{4-t^2,T-t^2}'(Z_{t^2})|  
 , \]
where $Z_{t^2} = Z_{t^2}(i\delta)$.  The result then
follows from \eqref{good.4}.
\end{proof}

\begin{lemma}  
For all subpower functions $\psi,
\phi$ there
is a subpower function $\overline \psi$ such that 
if $T$ is $\psi$-good at $\delta$, then
the following holds for $1 \leq t \leq 2$.
Suppose $w = x+iy$ with
\[    \delta \leq y \leq 1, \;\;\;\;\;
      (x/y)^2 + 1 \leq \phi(1/y).\]
Then,
\[   y^{-\beta} \, \overline \psi_\delta(1/y)^{-1} 
\leq   |h_{t^2}'(w)| \leq y^{-\beta} \, \overline \psi_\delta(1/y).\]
In particular,
\begin{equation}  \label{apr7.2}
  y^{-\beta} \, \overline \psi_\delta(1/y)^{-1} 
\leq   |h_{t^2}'(iy)| \leq y^{-\beta} \, \overline \psi_\delta(1/y).
\end{equation}
\end{lemma}

\begin{proof} We will do the case $t=2$; the argument
is similar for $1 \leq t \leq 2$.
 We let $\overline \psi$ denote a subpower
function in this proof, but its value may change
from line to line.  Since $(x/y)^2 + 1 \leq \phi
 (1/y)$,
we can see from Proposition \ref{prop34}
that
\[   \overline \psi(1/y)^{-1}  \, |h_{t^2}'(iy)| \leq
   |h_{t^2}'(w)| \leq \overline \psi(1/y) \, |h_{t^2}'(iy)|, \]
so we may assume $x = 0$.  Let $s = y$.  Using
the Loewner equation \eqref{reverse.1}, we can see that
there is a $c< \infty$ such that
\[       y \leq \Im h_{s^2}(iy) \leq c y , \]
\[       |\Re  h_{s^2}(iy)| \leq c y , \]
\[    |h_{s^2}'(iy)| \leq c . \]
The last estimate and \eqref{apr7.1} imply that
\[ |h_4'(y)| \asymp |{h}_{4-s^2,T-s^2}(Z_{s^2}(iy))|. \]
Using \eqref{good.5}, we see that
\[       \left| \Re Z_{s^2}(iy)\right|
    \leq y \, [c +\psi(1/y)]  . \]
By the previous result,
\[  \overline \psi(1/y)^{-1} \,t^{-\beta}\leq
 | h_{4- s^2,T-s^2}' (Z_{s^2}(iy))|
  \leq \overline \psi(1/y) \,t^{-\beta}
 \]
\end{proof}

\begin{definition}
If $n$ is a positive integer and $j=1,\ldots,n^2$,
we say that  $(j,n)$ is
$\psi$-good if $T= 1 + (j-1)n^{-2}$ is $\psi$-good
at $n^{-1}$.  We let  $E_{j,n}$ denote   the
event ``$(j,n)$ is
$\psi$-good '' 
and $E(j,n)$ denotes  the indicator function of 
$E_{j,n}$ .
\end{definition}

It is important to note that on the event $E_{j,n}$,
\eqref{apr7.2} implies that
\eqref{sept17.1}, the corresponding estimate for $|\hat{f}'|$, holds, with perhaps a different choice
of subpower function $\psi$. 
The main estimate for the lower bound is the following.

\begin{proposition}  \label{nov9.10}
If $r < r_c$,
here exists a subpower function $\psi$ and
$c > 0$ such that 
for all $n$ and all $j=1,2,\ldots,n^2$, 
\begin{equation}  \label{sept29.1}
  \E\left[|\hat f_{j,n}'(i/n)|^\lambda
  \, E_{j,n}\right]  \geq c \, n^{-\zeta} . 
\end{equation}
\end{proposition}

\begin{remark}  For fixed $n$, the expectation in
\eqref{sept29.1} is the
same for all $j$. 
\end{remark}

We will not include a proof of Proposition \ref{nov9.10} because
it has essentially appeared in  \cite[Theorem 10.8]{Law1}; 
see also \cite[Lemma 4.4]{JL}. We
point out that the assumption $r < r_c$ is crucial for the result.
The proof in  \cite{Law1} and \cite{JL} uses
a careful analysis of a relatively simple
one-dimensional diffusion.

\subsection{Correlations}

In this subsection we fix a subpower function
$\psi$ such that Proposition \ref{nov9.10}
holds.   If $n$ is a positive
integer, we write $j,k$ for positive integers satisfying
$1 \leq j < k \leq n^2$.   We will consider
$E_{j,n} \cap E_{k,n}$ with indicator function
$E(j,n) \, E(k,n).$   If $j,k,n$ are fixed we write
\[     S = 1 + \frac{j-1}{n^2}, \;\;\;\;\;
     T = 1 + \frac{k-1}{n^2} , \]
\[  \tilde h_{t} = h_{t,S}, \;\;\;\;  h_t =
     h_{t,T},\]
     and recall that this means that $\tilde{h}_t,  h_t$ are solutions to the reverse-time Loewner equation with $V_{S-t}-V_S$ and $V_{T-t}-V_T$ as driving functions, respectively.


\begin{proposition} \label{nov9.10.1} There is a subpower function
$\phi$ such that for all $1 \leq j < k \leq n^2$
\[  \E\left[|\hat f_T'( i/n)|^\lambda \, |\hat f_{S}'
   (i/n)|^\lambda \, E(j,n) \, E(k,n) \right]
    \leq n^{-2 \zeta} \, \left(\frac{n^2}{k-j} \right)
        ^{\frac{\lambda \beta+\zeta}{2}} \, \phi\left(
    \frac{n^2}{{k-j}}\right) . \]
Moreover, on the event $E_{j,n} \cap E_{k,n}$,
\[   |\hat f_T(i/n) - \hat f_S(i/n)| \geq 
   \left(\frac{{k-j}}{n^2}\right)^{\frac{1-\beta}
 2} \, 
  \phi\left(\frac{n^2}{k-j}\right) \]
\end{proposition}

\begin{proof}  We write $\phi$ for a subpower function
but we let its value vary from line to line; in the end
we choose the maximum of all the subpower functions
mentioned.
Recall that $\hat f_S'(i/n)  = \tilde h_S'(i/n),
\hat f_T'(i/n) = h_T'(i/n)$ and
\[  \hat f_T(i/n) - \hat f_S(i/n) =
     \tilde h_S(Z_{T-S}(i/n)) - \tilde h_S
   (i/n). \]
The second assertion of the proposition
follows immediately from
Lemma \ref{nov9.lemma4}, so we need only show the first.

Since $T$ is $\psi$-good at $1/n$, we know 
from \eqref{good.4} that
\[   |h_{T}'(i/n)| \leq |h_{T-S}'(i/n)|
  \, \left(\frac{n^2}{
    k-j}\right)^{\beta/2}  \, \phi
  \left(\frac{n^2}{
    k-j}\right). \]
Therefore,
\[ \E\left[|\hat f_T'( i/n)|^\lambda \, |\hat f_{S}'
   (i/n)|^\lambda \, E(j,n) \, E(k,n) \right]
\leq \hspace{1.5in} \]
\[ \hspace{.2in} \left(\frac{n^2}{
    k-j}\right)^{\lambda \beta/2}  \, \phi
  \left(\frac{n^2}{
    k-j}\right)\,
\E\left[|h_{T-S}'( i/n)|^\lambda \, |\tilde h_{S}'
   (i/n)|^\lambda \, E(j,n) \, E(k,n) \right]
.\]
Note that
$ |h_{T-S}'( i/n)|^\lambda \, E(k,n)$ and
$ |\tilde h_{S}'
   (i/n)|^\lambda \, E(j,n)$ are independent random variables.
Therefore,
\[ \E\left[|h_{T-S}'( i/n)|^\lambda \, |\tilde h_{S}'
   (i/n)|^\lambda \, E(j,n) \, E(k,n) \right]
 \leq \hspace{1.5in}\]
\[ \hspace{.8in} 
 \E\left[ |\tilde h_{S}'
   (i/n)|^\lambda \, E(j,n)\right] \; 
  \E\left[|h_{T-S}'(i/n)|^\lambda \,   E(k,n) \right]
.\]
We now apply Theorem \ref{nov9.prop20} to see that
the right hand side above is bounded above by
\[      n^{-\zeta} \, \left({k-j}\right)
   ^{-\zeta/2} \, \phi\left(\frac{n^2}{k-j}\right)
  =  n^{-2\zeta} \, \left(\frac{n^2}{k-j}\right)
   ^{\zeta/2} \, \phi\left(\frac{n^2}{k-j}\right), \]
 and this concludes the proof.
 \end{proof}

\section{Proof of Theorem \ref{newtheorem}}  \label{newsec}

In this section we will use the \emph{forward} Loewner flow to prove
Theorem \ref{newtheorem}, which we restate for the convenience of the reader.
\begin{theorem*}  
If $0< \kappa < 8$ and $1/2 \le \alpha \le \alpha_*$, then with probability one
there exists a set $V$ such that  $\hdim[\gamma(V)] \le F_{\rm{tip}}(\alpha)$
and for $t \not\in V$, $\gamma(t) \in \Half$, 
\begin{equation*} 
    \tilde \mu(t,2^{-n}) \lap  2^{-n\alpha}, 
  \;\;\; n \rightarrow \infty.
\end{equation*}  
\end{theorem*}

 Throughout
we will fix $\kappa = 2/a < 8$.  We will write $u$ rather than
$\alpha$ (to avoid having both $\alpha$ and $a$ in formulas). 
To prove the theorem
it suffices to show that for every bounded domain $D \subset
\Half$ bounded away from the real line, there is a set $V_D$
with $\hdim[\gamma(V_D)] \le F_{\rm{tip}}(u)$ and  such
that \eqref{aug13.5} holds for $t \not \in V_D$ with $\gamma(t)
 \in D$.  We fix such a $D$ and allow constants to depend
on $D$. The basic strategy is typical for establishing upper
bounds for multifractal spectra.  We estimate a particular
moment of $|g_\tau'(z)|$ for an appropriate stopping time,
use Chebyshev's inequality to get an estimate on probabilities,
and use this estimate to bound the dimension of a well chosen
covering.

 We warn the reader again that
 some of the notation in this section
is not consistent with that in other sections.

 We parametrize $\SLE_\kappa$ so
that the conformal maps $g_t$ satisfy
\begin{equation}  \label{forloew}
   \p_t g_t(z) = \frac{a}{g_t(z) - U_t},\;\;\;\;
  g_0(z) = z,
\end{equation}
where $U_t = -B_t$ is a standard Brownian motion.  This
is valid for $z \in \C \setminus \{0\}$ up to time
$T_z \in (0,\infty]$.  We let $H_t$ be the unbounded
component of $\Half \setminus \gamma(0,t]$.

\subsection{Preliminaries}

Let
\[  Z_t = Z_t(z) = X_t + iY_t = g_t(z) - U_t.\]
If $z \in \Half$,   let
\[  \deriv_t = |g_t'(z)|, \;\;\;\;
     \Upsilon_t =  \frac{Y_t}{|g_t'(z)|} 
,\;\;\;\;\Theta_t = \arg Z_t,\;\;\;\;
   S_t = \sin \Theta_t.\]
$\Upsilon_t$ equals $1/2$ times the conformal radius
of $H_t$ about $z$ (or we can
think of it as the conformal radius normalized so that
the conformal radius of $\Half$ about $i$ equals
1).   The Koebe-1/4 theorem implies that 
\begin{equation}  \label{koebenew}
   \frac { \Upsilon_t}{2} \leq
  \dist(z, \p H_t) \leq 2\, \Upsilon_t.
\end{equation}
Straightforward computations using 
\eqref{forloew} show that   for $z \in \Half$,
\[  \p_t \Delta_t =   \Delta_t \, \frac{a(Y_t^2
             - X_t^2)}{|Z_t|^4}, \;\;\;\;
   \p_t \Upsilon_t = -\Upsilon_t \, \frac{2a Y_t^2}
         {|Z_t|^4}
. \]
There exists $0 < c_1 < c_2 < \infty$ such that for $z \in D$,
\begin{equation} \label{sep30.1}
   c_1 \leq \Upsilon_0 \leq c_2, \;\;\;\;
     c_1 \leq S_0 \leq 1. 
\end{equation}

Let for $u \ge 1/2$
\begin{equation}  \label{aug12.1}
  r = r(u)  =   \frac 12   -2a
  - \frac 1{4u-2},
\end{equation}
\begin{equation}  \label{1}
  \lambda = \lambda_r = \frac{r^2}{2a} + r \, \left(
   1 - \frac 1{2a}\right), \;\;\;\;\xi =
\xi_r = \frac{r^2}{4a} = \frac \lambda 2 - \frac r2
  \, \left(
   1 - \frac 1{2a}\right).
\end{equation}
Note that $r$ increases with $u$. Define
\[
\hat{u}_c=\frac12 +\frac{1}{8a-2}.
\] 
Note that $\hat{u}_c < \alpha_*=2a/(4a-1)$.
If $u < \hat{u}_c$, since $a > 1/4$, 
\[    r < r(\hat{u}_c) = 1-4a < \min\left\{\frac
 12 -2a,2-3a\right\},\]\[
       r< 0, \;\;\;\; r + \lambda > 0. \] 
The following is a straightforward It\^o's formula
calculation that we omit.

\begin{proposition} \label{aug12.prop1}
 Suppose $r \in \R$ and $\lambda,
\xi$ are as in \eqref{1}.  
If $z \in \Half$, and 
\[            M_t =  M_t(z) =
|Z_t|^r \, Y_t^\xi \, \deriv_t^\lambda
   = S_t^{-r} \, \Upsilon_t^{\xi + r} \, \deriv_t^{\lambda
  + r}, \]
then $M_t$ is a local martingale satisfying
\[     dM_t = M_t \, \frac{r \, X_t}{|Z_t|^2} \, dB_t.\]
\end{proposition}

Let $\dyad_n$ denote the of dyadic rationals
in $\C$
\[         z=  \frac{j}{2^n} + i \frac{k}{2^n}, \;\;\;\;
  j,k \in \Z.\]
Note that if
 $w \in \C$,  then there exists
$z \in \dyad_n$ with $|z-w| \leq 2^{-n}$ and
hence $\ball(w,2^{-n}) \subset \ball(z,2^{-n+1}).$

\subsection{Basic strategy}

Let
\[  \tau_{n,z} = \inf\left\{s: \Upsilon_s(z)  \leq 2^{-n+3}
  \right\}. \]
We will only consider $n$ sufficiently large  so that
$2^{-n+4} \leq c_1$ where $c_1$ is the constant in
\eqref{sep30.1}.  Note that $\Prob\{\tau_{n,z} = \infty\}
>0$. If  $\tau_{n,z} < \infty
 $, \eqref{koebenew} implies
\[      2^{-n+2} \leq \dist(z,\p H_{\tau_{n,z}}) =
   \dist(z,\gamma(0,\tau_{n,z}]) \leq 2^{-n+4}. \]
In particular, if $|w-z| \leq 2^{-n}, $
 $\dist(w,\p H_{\tau_{n,z}}) \geq 2^{-n+1}$.
 
Recall that we defined the normalized harmonic measure
\[
\hm_t(V)=\lim_{y \to \infty} y \, \hm(iy, V, H_t)
\]
and $\tilde{\mu}(t, \epsilon)=\hm_t(\overline{\mathcal{B}}(\gamma(t), \epsilon))$. Similarly we define
\[
\hat{\hm}_t(V):=\lim_{y \to \infty} y \, \hm(iy, \partial V, H_t \setminus V),
\]
and note that $\hm_t(V) \le \hat{\hm}_t(V)$. Set $ \hm_{n,z} = \hm_{\tau_{n,z}}$ and similarly for $\hat{\hm}$.  If $|z - \gamma(t)| \leq
2^{-n}$, then $\tau_{n,z} \leq t$ and hence by monotonicity
of harmonic measure, 
\[      \hat{\hm}_{n,z} \left[\ball(z,2^{-n+1})\right]
           \geq {\hm}_t\left[\ball(z,2^{-n+1})\right]
  \geq \tilde \mu(t,2^{-n}) .\]

Let $\dyad_n(D)$ denote the set of $z \in \dyad$ such
that $\dist(z,D) \leq 2^{-n}$ and 
\[   A^{u}_m =  A^u_m(D) =
\bigcup_{n=m}^\infty  \bigcup_z \ball(z,2^{-n+1}) , \]
where the inner union is over all $z \in \dyad_n(D)$ satisfying
\begin{equation}  \label{aug11.2}
         \hat{\hm}_{n,z} \left[\ball(z,2^{-n+1})\right]
   \geq 2^{-nu}. 
\end{equation}
Then if $ \gamma(t)  \in D \setminus
A^u_m$, for all $n$ sufficiently large
\[   {\tilde{\mu}}(t, 2^{-n})  \leq 2^{-nu}.\]
Hence for each $m$, $A_m^u$ is a cover of  
$D\cap V_u$ where $V_u$ is  the set 
of $\gamma(t)$ that do not satisfy \eqref{aug13.5}.
   Let $N_n = N_{n,u}(D)$ be
the cardinality of the set of $z \in \dyad_n(D)$ satisfying \eqref{aug11.2}.  Then
for all $s$, 
\[  \haus^s\left[D  \cap
V_u\right]
 \leq \lim_{m \rightarrow \infty} 
\haus^s\left( A^{u}_m \right)
            \leq c\lim_{m \rightarrow \infty}
\sum_{n=m}^\infty N_n \, 2^{-ns}. \]
The following proposition follows immediately.

\begin{proposition}  Suppose   $u,s > 0$
and
\[              N_{n,u}(D) \lap 2^{ns}, \;\;\;\;
  n \rightarrow \infty. \]
Then
\[         \hdim\left[D \cap V_u\right]
   \leq  s. \]
\end{proposition}

\begin{proof}  The argument above shows that for all $s' >s$,
 $\haus^{s'}\left[D \cap V_u \right]=0$.
\end{proof}

In order to show that with probability one for $\SLE$ that
\[         \hdim\left[ D \cap V_u\right]
   \leq  s. \]
it suffices to show that
\[    \E\left[  N_{n,u}(D) \right] \lap 2^{ns}, \;\;\;\;
  n \rightarrow \infty. \]
Indeed, this relation and the Borel-Cantelli lemma imply
that with probability one for all $s' > s$,
$N_{n,u}(D) \leq 2^{ns'}.$  Note that
\[  \E\left[  N_{n,u}(D) \right]  \leq c_D\,  2^{2n}
  \, \sup_{\dist(z,D) \leq 2^{-n}}  \Prob
 \left\{\tau_{n,z} < \infty; 
         \, \hat{\hm}_{n,z}[\ball(z,2^{-n+1}) ] \geq 2^{-nu} \right\}. \]

Notice that conformal invariance of harmonic measure and distortion
estimates imply that on the event $\tau_{n,z} < \infty$,
\[      \hat{\hm}_{n,z} \left[\ball(z,2^{-n+1})\right]
         \asymp   2^{-n} \, |g_{\tau_{n,z}}'(z)|   
 . \]
Indeed, $g_{\tau_{n,z}}(\ball(z,2^{-n+1}))$ is a connected
set whose diameter is
comparable to $2^{-n} \, |g_{\tau_{n,z}}'(z)|$ and whose  distance
from the real axis is 
comparable to $2^{-n} \, |g_{\tau_{n,z}}'(z)|$.
Hence, there exists $c < \infty$ such that 
 \[  \E\left[  N_{n,u}(D) \right]  \leq c_D\,  2^{2n}
  \, \sup_{\dist(z,D) \leq 2^{-n}}  \Prob
 \left\{\tau_{n,z} < \infty; 
         |g_{\tau_{n,z}}'(z)| \geq c \,  2^{-n(u-1)} \right\}. \]

In the remainder of this section we will show that
there exists $c = c_D < \infty$ such that for all $n$
sufficiently large and all $z$ with $\dist(z,D) \leq 2^{-n}$,
\begin{equation}  \label{aug13.6}
  \Prob
 \left\{\tau_{n,z} < \infty; 
         |g_{\tau_{n,z}}'(z)| \geq c \,  2^{-n(u -
 1)} \right\}
 \leq 2^{-n\rho(u)}, 
\end{equation}
where
\[    \rho(u) = \left[\frac{1}{8a} + 2a - 1
 \right] \, ( u-\frac{1}{2} ) + \left[\frac 12 - \frac 1{8a} 
\right] + \frac{1}{32a(u-\frac{1}{2})}. \]
Then from the above arguments we know
that with probability one,
\[         \hdim\left[D \cap V_u\right]
   \leq  2 - \rho(u) = F_{\rm tip}(u). \]
The second equality is a straight-forward calculation.  The remainder
of this section  is devoted to establishing
\eqref{aug13.6}.

\subsection{Weighting by the martingale}

The local martingale $M_t$ is not a martingale because it ``blows
up'' on the event of measure zero that $z$ is
on the path  $\gamma(0,\infty)$.
However, if we choose stopping times $\tau$
such as $\tau_{n,z}$ which
prevent the path from getting too close to $z$, then the
stopped process $M_{t \wedge \tau}$ is a martingale.
 Let $\Prob^*,\E^*$
denote  probabilities and
expectations with respect to the measure obtained by weighting by (the stopped martingale) $M$.
The Girsanov theorem implies that
\[         dB_t = \frac{r \,X_t}{|Z_t|^2} \, dt + dW_t, \;\;\;\;
 0 \leq t < \tau,  \]
where $W_t$ is a standard Brownian motion with respect
to the measure $\Prob^*$. 
In particular,
\[   d \Theta_t = \frac{(1-2a-r) \, X_t\, Y_t}{|Z_t|^4}
  \, dt  - \frac{Y_t}{|Z_t|^2} \, dW_t. \]

It is useful to use a ``radial'' parametrization
  $\sigma(t)$.  We   
write $\hat Z_t = Z_{\sigma(t)},$ $\hat{X}_t $ $= X_{\sigma(t)},$
etc. 
The radial parametrization is defined by
\[        \hat \Upsilon_t :=
     \Upsilon_{\sigma(t)} = e^{-2at}. \]
Note that
\[  -2a \, \hat \Upsilon_t = \p_t \hat \Upsilon_t
        = -2a \, \frac{\hat Y_t^2}{|\hat Z_t|^4} \, 
  \p_t \sigma(t) , \]
which implies
\[           \p_t \sigma(t) = \frac{|\hat Z_t|^4}{\hat{Y}_t^2}.\]
Note also that
\[  d \hat \Theta_t = (1-2a) \, \cot \hat \Theta_t \, dt +
           d \hat B_t, \]
and the local martingale $M_t$ satisfies
\[          d \hat
 M_t = - r \, \hat M_t \,  \cot \hat \Theta_t \, d
  \hat B_t.\] Moreover, we have that
\begin{equation} \label{sep30.2}
   d \hat \Theta_t = (1-2a-r) \, \cot \hat \Theta_t \, dt
            + d \hat W_t.   
\end{equation}
In the above, $\hat B_t$ and $\hat W_t$ are standard Brownian motions
with respect to $\Prob$ and $\Prob^*$ respectively.
 Since $1-2a-r > 1/2$, we compare with a Bessel process to see that in the measure $\Prob^*$,
 $\hat \Theta_t$ never reaches $\{0,\pi\}$, see \cite[Chapter 1]{LConv}. It follows that $\hat M_t$ is actually a martingale. Also,
the invariant probability density for the SDE
\eqref{sep30.2} equals
\[    f(\theta) = c\, \sin^{2(1-2a-r)} \theta. \]
Since $r < 1-4a < 3-4a$ it follows that $\sin^r$ is integrable with respect to
$f(\theta) \, d\theta$.  The important fact for us,
is that there is a $c$ such that  if
$\hat \Theta_t$ satisfies \eqref{sep30.2} with
  $\sin \hat \Theta_0 \geq c_1$, then for all $t > 0$,
\begin{equation}  \label{aug13.7}
     \E^*[\hat S_t^r ]  \leq c. 
\end{equation}

Let 
\[  \tau_s = \inf
  \left \{t: \Upsilon_t = e^{-2as} \right\}.\]
For $r < \frac 12 - 2a, $ we have for all $s$,
\[  \Prob^*\{\tau_s < \infty \} = 1.\]
Then, using \eqref{aug13.7}, 
\begin{eqnarray*}
  \E\left[|g_{\tau_s}'(z)|^{\lambda+r};
           \tau_s < \infty \right]& 
  = & \E\left[M_{\tau_s} \, S_{\tau_s}^{r}
  \, \Upsilon_{\tau_s}^{-\xi-r}; \tau_s < \infty\right] \\
 & = & e^{2as(\xi+r)} \, 
 \E\left[M_{\tau_s} \, S_{\tau_s}^{r}
  ; \tau_s < \infty\right] \\
& = & e^{2as(\xi+r)} \, 
M_0(s) \, \E^*\left[\hat S_{s}^r\right]\\
& \leq & c \, e^{2as(\xi+r)}. 
\end{eqnarray*}

Since $\lambda + r > 0$,  if $\epsilon = e^{-2as}$,
\[  \Prob\left\{\tau_s < \infty; 
|g_{\tau_s}'(z)| \geq \epsilon^{u- 1} \right\}
   \leq  \epsilon^{-(u- 1)(\lambda+r)} 
 \, \E[|g_{\tau_s}'(z)|^{\lambda+r}; \tau_s < \infty]
  \leq  c\, \epsilon^{\rho(u) }\]
where
\[ \rho(u) = -(u- 1)(\lambda+r) - (r + \xi). \]
Doing the algebra, we get
\[    \rho(u) = \left[\frac{1}{8a} + 2a - 1
 \right] \, (u-\frac{1}{2}) + \left[\frac 12 - \frac 1{8a} 
\right] + \frac1{32a(u-\frac{1}{2})}. \]
 This proves \eqref{aug13.6} which concludes the proof of Theorem~\ref{newtheorem}. \flushright{$\square$}

\end{document}